\documentclass[11pt]{article}
\usepackage{amsmath, upgreek}
\usepackage{amssymb}
\usepackage{color}
\usepackage{amscd}
\usepackage{xspace}
\usepackage{verbatim}
\usepackage{graphicx}
\renewcommand{\theequation}{\thesection.\arabic{equation}
}
\newcommand{\mysection}[1]{
\section{#1}\setcounter{equation}{0}}
\title{\bf Boundary singularities of  semilinear elliptic \\ equations   
with Leray-Hardy potential}
\author{{\bf Huyuan Chen\footnote{\noindent Department of Mathematics, Jiangxi Normal University,
Nanchang 330022, China. E-mail: chenhuyuan@yeah.net}} \\[4mm]
 {\bf Laurent V\'eron \footnote{\noindent
Laboratoire de Math\'{e}matiques et Physique Th\'{e}orique, Universit\'e de Tours, 37200 Tours, France. E-mail: veronl@univ-tours.fr}}
}

\date{}
\begin{document}
 \maketitle


\newcommand{\txt}[1]{\;\text{ #1 }\;}
\newcommand{\tbf}{\textbf}
\newcommand{\tit}{\textit}
\newcommand{\tsc}{\textsc}
\newcommand{\trm}{\textrm}
\newcommand{\mbf}{\mathbf}
\newcommand{\mrm}{\mathrm}
\newcommand{\bsym}{\boldsymbol}
\newcommand{\scs}{\scriptstyle}
\newcommand{\sss}{\scriptscriptstyle}
\newcommand{\txts}{\textstyle}
\newcommand{\dsps}{\displaystyle}
\newcommand{\fnz}{\footnotesize}
\newcommand{\scz}{\scriptsize}
\newcommand{\be}{\begin{equation}}
\newcommand{\bel}[1]{\begin{equation}\label{#1}}
\newcommand{\ee}{\end{equation}}
\newcommand{\eqnl}[2]{\begin{equation}\label{#1}{#2}\end{equation}}
\newcommand{\barr}{\begin{eqnarray}}
\newcommand{\earr}{\end{eqnarray}}
\newcommand{\bars}{\begin{eqnarray*}}
\newcommand{\ears}{\end{eqnarray*}}
\newcommand{\nnu}{\nonumber \\}
\newtheorem{subn}{\name}
\renewcommand{\thesubn}{}
\newcommand{\bsn}[1]{\def\name{#1}\begin{subn}}
\newcommand{\esn}{\end{subn}}
\newtheorem{sub}{\name}[section]
\newcommand{\dn}[1]{\def\name{#1}}   
\newcommand{\bs}{\begin{sub}}
\newcommand{\es}{\end{sub}}
\newcommand{\bsl}[1]{\begin{sub}\label{#1}}
\newcommand{\bth}[1]{\def\name{Theorem}
\begin{sub}\label{t:#1}}
\newcommand{\blemma}[1]{\def\name{Lemma}
\begin{sub}\label{l:#1}}
\newcommand{\bcor}[1]{\def\name{Corollary}
\begin{sub}\label{c:#1}}
\newcommand{\bdef}[1]{\def\name{Definition}
\begin{sub}\label{d:#1}}
\newcommand{\bprop}[1]{\def\name{Proposition}
\begin{sub}\label{p:#1}}

\newcommand{\aand}{\quad\mbox{and}\quad}
\newcommand{\M}{{\cal M}}
\newcommand{\A}{{\cal A}}
\newcommand{\B}{{\cal B}}
\newcommand{\I}{{\cal I}}
\newcommand{\J}{{\cal J}}
\newcommand{\D}{\displaystyle}
\newcommand{\RR}{ I\!\!R}
\newcommand{\C}{\mathbb{C}}
\newcommand{\R}{\mathbb{R}}
\newcommand{\Z}{\mathbb{Z}}
\newcommand{\N}{\mathbb{N}}
\newcommand{\T}{{\rm T}^n}
\newcommand{\cuad}{{\sqcap\kern-.68em\sqcup}}
\newcommand{\abs}[1]{\mid #1 \mid}
\newcommand{\norm}[1]{\|#1\|}
\newcommand{\equ}[1]{(\ref{#1})}
\newcommand\rn{\mathbb{R}^N}
\renewcommand{\theequation}{\thesection.\arabic{equation}}
\newtheorem{definition}{Definition}[section]
\newtheorem{theorem}{Theorem}[section]
\newtheorem{proposition}{Proposition}[section]
\newtheorem{example}{Example}[section]
\newtheorem{proof}{proof}[section]
\newtheorem{lemma}{Lemma}[section]
\newtheorem{corollary}{Corollary}[section]
\newtheorem{remark}{Remark}[section]
\newcommand{\bremark}{\begin{remark} \em}
\newcommand{\eremark}{\end{remark} }
\newtheorem{claim}{Claim}


\newcommand{\rth}[1]{Theorem~\ref{t:#1}}
\newcommand{\rlemma}[1]{Lemma~\ref{l:#1}}
\newcommand{\rcor}[1]{Corollary~\ref{c:#1}}
\newcommand{\rdef}[1]{Definition~\ref{d:#1}}
\newcommand{\rprop}[1]{Proposition~\ref{p:#1}}
\newcommand{\BA}{\begin{array}}
\newcommand{\EA}{\end{array}}
\newcommand{\BAN}{\renewcommand{\arraystretch}{1.2}
\setlength{\arraycolsep}{2pt}\begin{array}}
\newcommand{\BAV}[2]{\renewcommand{\arraystretch}{#1}
\setlength{\arraycolsep}{#2}\begin{array}}
\newcommand{\BSA}{\begin{subarray}}
\newcommand{\ESA}{\end{subarray}}
\newcommand{\BAL}{\begin{aligned}}
\newcommand{\EAL}{\end{aligned}}
\newcommand{\BALG}{\begin{alignat}}
\newcommand{\EALG}{\end{alignat}}
\newcommand{\BALGN}{\begin{alignat*}}
\newcommand{\EALGN}{\end{alignat*}}
\newcommand{\note}[1]{\textit{#1.}\hspace{2mm}}
\newcommand{\Proof}{\note{Proof}}
\newcommand{\qeda}{\hspace{10mm}\hfill $\square$}
\newcommand{\qed}{\\
${}$ \hfill $\square$}
\newcommand{\Remark}{\note{Remark}}
\newcommand{\modin}{$\,$\\[-4mm] \indent}
\newcommand{\forevery}{\quad \forall}
\newcommand{\set}[1]{\{#1\}}
\newcommand{\setdef}[2]{\{\,#1:\,#2\,\}}
\newcommand{\setm}[2]{\{\,#1\mid #2\,\}}
\newcommand{\mt}{\mapsto}
\newcommand{\lra}{\longrightarrow}
\newcommand{\lla}{\longleftarrow}
\newcommand{\llra}{\longleftrightarrow}
\newcommand{\Lra}{\Longrightarrow}
\newcommand{\Lla}{\Longleftarrow}
\newcommand{\Llra}{\Longleftrightarrow}
\newcommand{\warrow}{\rightharpoonup}
\newcommand{
\paran}[1]{\left (#1 \right )}
\newcommand{\sqbr}[1]{\left [#1 \right ]}
\newcommand{\curlybr}[1]{\left \{#1 \right \}}
\newcommand{
\paranb}[1]{\big (#1 \big )}
\newcommand{\lsqbrb}[1]{\big [#1 \big ]}
\newcommand{\lcurlybrb}[1]{\big \{#1 \big \}}
\newcommand{\absb}[1]{\big |#1\big |}
\newcommand{\normb}[1]{\big \|#1\big \|}
\newcommand{
\paranB}[1]{\Big (#1 \Big )}
\newcommand{\absB}[1]{\Big |#1\Big |}
\newcommand{\normB}[1]{\Big \|#1\Big \|}
\newcommand{\produal}[1]{\langle #1 \rangle}

\newcommand{\thkl}{\rule[-.5mm]{.3mm}{3mm}}
\newcommand{\thknorm}[1]{\thkl #1 \thkl\,}
\newcommand{\trinorm}[1]{|\!|\!| #1 |\!|\!|\,}
\newcommand{\bang}[1]{\langle #1 \rangle}
\def\angb<#1>{\langle #1 \rangle}
\newcommand{\vstrut}[1]{\rule{0mm}{#1}}
\newcommand{\rec}[1]{\frac{1}{#1}}
\newcommand{\opname}[1]{\mbox{\rm #1}\,}
\newcommand{\supp}{\opname{supp}}
\newcommand{\dist}{\opname{dist}}
\newcommand{\myfrac}[2]{{\displaystyle \frac{#1}{#2} }}
\newcommand{\myint}[2]{{\displaystyle \int_{#1}^{#2}}}
\newcommand{\mysum}[2]{{\displaystyle \sum_{#1}^{#2}}}
\newcommand {\dint}{{\displaystyle \myint\!\!\myint}}
\newcommand{\q}{\quad}
\newcommand{\qq}{\qquad}
\newcommand{\hsp}[1]{\hspace{#1mm}}
\newcommand{\vsp}[1]{\vspace{#1mm}}
\newcommand{\ity}{\infty}
\newcommand{\prt}{\partial}
\newcommand{\sms}{\setminus}
\newcommand{\ems}{\emptyset}
\newcommand{\ti}{\times}
\newcommand{\pr}{^\prime}
\newcommand{\ppr}{^{\prime\prime}}
\newcommand{\tl}{\tilde}
\newcommand{\sbs}{\subset}
\newcommand{\sbeq}{\subseteq}
\newcommand{\nind}{\noindent}
\newcommand{\ind}{\indent}
\newcommand{\ovl}{\overline}
\newcommand{\unl}{\underline}
\newcommand{\nin}{\not\in}
\newcommand{\pfrac}[2]{\genfrac{(}{)}{}{}{#1}{#2}}

\def\ga{\alpha}     \def\gb{\beta}       \def\gg{\gamma}
\def\gc{\chi}       \def\gd{\delta}      \def\ge{\epsilon}
\def\gth{\theta}                         \def\vge{\varepsilon}
\def\gf{\phi}       \def\vgf{\varphi}    \def\gh{\eta}
\def\gi{\iota}      \def\gk{\kappa}      \def\gl{\lambda}
\def\gm{\mu}        \def\gn{\nu}         \def\gp{\pi}
\def\vgp{\varpi}    \def\gr{\rho}        \def\vgr{\varrho}
\def\gs{\sigma}     \def\vgs{\varsigma}  \def\gt{\tau}
\def\gu{\upsilon}   \def\gv{\vartheta}   \def\gw{\omega}
\def\gx{\xi}        \def\gy{\psi}        \def\gz{\zeta}
\def\Gg{\Gamma}     \def\Gd{\Delta}      \def\Gf{\Phi}
\def\Gth{\Theta}
\def\Gl{\Lambda}    \def\Gs{\Sigma}      \def\Gp{\Pi}
\def\Gw{\Omega}     \def\Gx{\Xi}         \def\Gy{\Psi}

\def\CS{{\mathcal S}}   \def\CM{{\mathcal M}}   \def\CN{{\mathcal N}}
\def\CR{{\mathcal R}}   \def\CO{{\mathcal O}}   \def\CP{{\mathcal P}}
\def\CA{{\mathcal A}}   \def\CB{{\mathcal B}}   \def\CC{{\mathcal C}}
\def\CD{{\mathcal D}}   \def\CE{{\mathcal E}}   \def\CF{{\mathcal F}}
\def\CG{{\mathcal G}}   \def\CH{{\mathcal H}}   \def\CI{{\mathcal I}}
\def\CJ{{\mathcal J}}   \def\CK{{\mathcal K}}   \def\CL{{\mathcal L}}
\def\CT{{\mathcal T}}   \def\CU{{\mathcal U}}   \def\CV{{\mathcal V}}
\def\CZ{{\mathcal Z}}   \def\CX{{\mathcal X}}   \def\CY{{\mathcal Y}}
\def\CW{{\mathcal W}} \def\CQ{{\mathcal Q}}
\def\BBA {\mathbb A}   \def\BBb {\mathbb B}    \def\BBC {\mathbb C}
\def\BBD {\mathbb D}   \def\BBE {\mathbb E}    \def\BBF {\mathbb F}
\def\BBG {\mathbb G}   \def\BBH {\mathbb H}    \def\BBI {\mathbb I}
\def\BBJ {\mathbb J}   \def\BBK {\mathbb K}    \def\BBL {\mathbb L}
\def\BBM {\mathbb M}   \def\BBN {\mathbb N}    \def\BBO {\mathbb O}
\def\BBP {\mathbb P}   \def\BBR {\mathbb R}    \def\BBS {\mathbb S}
\def\BBT {\mathbb T}   \def\BBU {\mathbb U}    \def\BBV {\mathbb V}
\def\BBW {\mathbb W}   \def\BBX {\mathbb X}    \def\BBY {\mathbb Y}
\def\BBZ {\mathbb Z}

\def\GTA {\mathfrak A}   \def\GTB {\mathfrak B}    \def\GTC {\mathfrak C}
\def\GTD {\mathfrak D}   \def\GTE {\mathfrak E}    \def\GTF {\mathfrak F}
\def\GTG {\mathfrak G}   \def\GTH {\mathfrak H}    \def\GTI {\mathfrak I}
\def\GTJ {\mathfrak J}   \def\GTK {\mathfrak K}    \def\GTL {\mathfrak L}
\def\GTM {\mathfrak M}   \def\GTN {\mathfrak N}    \def\GTO {\mathfrak O}
\def\GTP {\mathfrak P}   \def\GTR {\mathfrak R}    \def\GTS {\mathfrak S}
\def\GTT {\mathfrak T}   \def\GTU {\mathfrak U}    \def\GTV {\mathfrak V}
\def\GTW {\mathfrak W}   \def\GTX {\mathfrak X}    \def\GTY {\mathfrak Y}
\def\GTZ {\mathfrak Z}   \def\GTQ {\mathfrak Q}

\font\Sym= msam10 
\def\SYM#1{\hbox{\Sym #1}}
\newcommand{\bdw}{\prt\Gw\xspace}
\date{}
\maketitle\medskip

\begin{abstract}
We study existence and uniqueness of solutions of ($E_1$) $-\Gd u+\tfrac {\gm}{|x|^2}u+g(u)=\gn$ in $\Gw$, $u=\gl$ on $\prt\Gw$,
where $\Gw\subset\BBR^N_+$ is a bounded smooth domain such that $0\in\prt\Gw$,  $\gm\geq-\frac{N^2}{4}$ is a constant, $g$ a continuous nondecreasing function satisfying some integral growth condition and 
$\gn$ and $\gl$ two Radon measures respectively in $\Gw$ and on $\prt\Gw$. We show that the situation differs considerably according the measure is concentrated at $0$ or not. When $g$ is a power we introduce a capacity framework which provides necessary and sufficient conditions for the solvability of problem ($E_1$).
\end{abstract}

\noindent
  \noindent {\small {\bf Key Words}:   Hardy Potential,  Radon Measure.  }\vspace{1mm}

\noindent {\small {\bf MSC2010}:  35B44, 35J75. }\tableofcontents
\vspace{1mm}
\hspace{.05in}
\medskip

\setcounter{equation}{0}
\section{Introduction}

If $\gm$ is a real number and $N\geq 2$,  the  Schr\"odinger operator $\CL_\gm$, defined in a domain (any connected open subset) $\Gw\subset\BBR^N$ by
\bel{1-A1}
\CL_\gm u:=-\Gd u+\frac{\gm}{|x|^2}u,
\ee
plays a fundamental role in analysis, because of Hardy's inequality, and in theoretical physics in connexion with uncertainty principle (see e.g. \cite {FS}, \cite {Fr}). 
When the singular point $0$ belongs to $\Gw$, there exists a critical value 
 \bel{1-A2}
\gm_0=-\left(\myfrac{N-2}{2}\right)^2.
\ee
It is wellknown that when $\gm\geq \gm_0$ the operator $\CL_\gm$ is positive because of Hardy's inequality
 \bel{1-A3}
\myint{\Gw}{}|\nabla \phi|^2dx+\gm_0\myint{\Gw}{}\myfrac{\phi^2}{|x|^2}dx\geq 0\quad\text{ for all }\gf\in C^\infty_0(\Gw).
\ee
Under the condition $\mu\geq \mu_0$, the study of semilinear problems associated to the operator $u\mapsto \CL_\gm u+g(u)$ has been initiated in \cite{GuVe}. When $g(r)\sim |r|^{p-1}r$ ($p>1$) the authors provided therein necessary and sufficient conditions on $p$  in order the no solution of equation
 \bel{1-A3'}\CL_\gm u+g(u)=0
\ee
in $\BBR^N\setminus\{0\}$ could have a singularity at $x=0$. When this condition is not satisfied they obtained a description of the possible behaviour of singular solutions of (\ref{1-A3'}) in the neighborhood of $0$. 
When $g$ is merely a continuous nondecreasing function they found a necessary and sufficient condition on $g$ expressed under an integral formulation ensuring the existence of a positive solution $u(x)$ of (\ref{1-A3'}) satisfying
 \bel{1-A3''}
 u(x)\sim a|x|^{-(\frac{N-2}{2}+\sqrt{\gm-\gm_0})} \quad\text{when }\;x\to 0,\; a>0.
\ee
Note that $x\mapsto |x|^{-(\frac{N-2}{2}+\sqrt{\gm+\gm_0})}$ is the solution of $\CL_\gm u=0$ in $\BBR^N\setminus\{0\}$ with the strongest singularity. Thanks to a notion of weak solutions of $\CL_\gm u=0$ combined with a  dual formulation of the equation introduced in \cite{CQZ},
we studied in \cite{ChVe} the equation
\bel{1-A4}\left.\BA {lll}
\CL_\gm u+ g(u)=\gn&\quad\text{in }\, \, \Gw\\[1mm]
\phantom{\CL_\gm + g(u)}
u=0&\quad\text{on }\, \prt\Gw,
\EA\right.\ee
in a bounded smooth domain $\Gw$, where $g$ is a continuous nondecreasing function and $\gn$ a Radon measure which support may contain $0$. In this framework, weak solutions to (\ref{1-A4}) with measures defined in a class of weighted measures are obtained  provided that $g$ satisfies some integrability condition. When this integrability condition is not satisfied by $g$, i.e. in  the {\it supercritical case}, not all measures in the above class are suitable for solving $(\ref{1-A4})$. In the particular case where $g(r)=|r|^{p-1}r$ ($p>1$), we proved that a weighted measure is suitable for solving (\ref{1-A4}) if it is absolutely continuous with respect to some specific Bessel capacity depending on $\gm$, $N$ and $p$. \smallskip

In this article we are interested in similar problems but in the configuration where the singular point $0$ of the Leray-Hardy potential lies on the boundary of the domain  $\Gw$. Our aim is to extend the approach developed in \cite{GuVe} and \cite{ChVe} for the problem with an internal singularity to the  study the following equation
\bel{1-A5}
 \left.\BA {lll}
\CL_\gm u+ g(u)=\gn&\quad\text{in }\, \Gw\\[1mm]
\phantom{\CL_\gm + g(u)}
u=\gl&\quad\text{on }\prt\Gw,
\EA\right.
\ee
where $\gn$ and $\gl$ are bounded Radon measures respectively on $\Gw$ and $\prt\Gw$.
When $\gm=0$ the first study is due to Gmira and V\'eron \cite{GmVe} who proved the existence and uniqueness of a {\it very weak} solution,  a denomination due to Brezis in an unpublished note. Such a solution $u$ is a function belonging to $L^1(\Gw)$ such that $\gr g(u)\in L^1(\Gw)$, where $\gr(x)=\dist (x,\prt\Gw)$, satisfying
 \bel{1-A6}
\BA {lll}
\myint{\Gw}{}\left(-u\Gd\gz+g(u)\gz\right) dx=-\myint{\prt\Gw}{}\myfrac{\prt\gz}{\prt{\bf n}}d\gn
\EA
\ee
for all $\gz\in C^{1}_c(\overline\Gw)$ such that $\Gd\gz\in L^\infty(\Gw)$ (in the sense of distributions in $\Gw$). A sufficient condition for the existence and uniqueness of a solution to (\ref{1-A5}) with $\gn=0$ is 
 \bel{1-A7}\BA {lll}
\myint{1}{\infty}\left(g(s)-g(-s)\right) s^{-\frac{2N}{N-1}}ds<\infty. 
\EA\ee
When $\gm\neq 0$,  a model domain  is 
$\Gw=\BBR^N_+:=\{x=(x',x_N)=(x_1,...,x_N)=x_N>0\}$. The operator $\CL_\gm$ is positive if   
\bel{1-A11}
\gm\geq\gm_1:=-\frac {N^2}{4},
\ee
and $\gm_1$ is the best constant of the Hardy inequality in $\BBR^{N}_+$ since there holds
\bel{1-A12*}
\myint{\BBR^{N}_+}{}|\nabla \phi|^2dx+\gm_1\myint{\BBR^{N}_+}{}\myfrac{\phi^2}{|x|^2}dx\geq 0\qquad\text{for all }\gf\in C^\infty_0(\BBR^{N}_+).
\ee 
If $\BBR^{N}_+$ is replaced by a bounded domain $\Gw$ satisfying the condition 
$$(\CC_1) \qquad\qquad\qquad 0\in \prt\Gw\text{ , }\;\Gw\subset\BBR^{N}_+\text{ and }\;\langle x,{\bf n}\rangle=O(|x|^2)\,\text{ for all }\;x\in\prt \Gw,\qquad\qquad\qquad\qquad\qquad\qquad$$
where ${\bf n}={\bf n}_x$ is the outward normal vector at $x$, this inequality is never achieved and there exists a remainder \cite{Caz1}: if we set 
$\displaystyle R_\Gw=\max_{z\in\Gw} |z|$, there holds
\bel{1-A12}
\myint{\Gw}{}|\nabla \phi|^2dx+\gm_1\myint{\Gw}{}\myfrac{\phi^2}{|x|^2}dx\geq \myfrac{1}{4}
\myint{\Gw}{}\myfrac{\phi^2}{|x|^2\ln^2(|x|R_\Gw^{-1})}dx
\qquad\text{for all }\gf\in C^\infty_0(\Gw).
\ee
Note that the last condition in $(\CC_1)$ holds if $\Gw$ is a $C^2$ domain. We put 
\bel{1-A14}
\ga_+:=\ga_+(\gm)=1-\frac{N}{2}+\sqrt{\gm+\frac{N^2}{4}}\quad\text{and }\;\ga_-:=\ga_-(\gm)=1-\frac{N}{2}-\sqrt{\gm+\frac{N^2}{4}}.
\ee
If $\Gw$ satisfies $(\CC_1)$  we define $\ell_\gm^\Gw$ by 
 \bel{1-A20}
 \BA {lll}
\ell_\gm^\Gw:=\min\left\{\myint{\Gw}{}\left(|\nabla v|^2+\myfrac{\gm}{|x|^2}v^2\right)dx:v\in C^{1}_c(\Gw),\myint{\Gw}{}v^2 dx=1
\right\}.\EA\ee
Then $\ell_\gm^\Gw>0$. If  $\gm\geq\gm_1$ this first eigenvalue of $\CL_\gm$ is achieved in the closure $ H_\gm(\Gw)$ of $C^{1}_c(\Gw)$  for the norm
\bel{1-A21}
 v\mapsto \norm v_{H_\mu(\Gw)}:=\sqrt{\myint{\Gw}{}\left(|\nabla v|^2+\myfrac{\gm}{|x|^2}v^2\right)dx}.
\ee
Note that $H_\gm(\Gw)=H^1_0(\Gw)$ if $\gm>\gm_1$, $H^1_0(\Gw)\subsetneqq  H_{\gm_1}(\Gw)$ and the imbedding of $H_{\gm_1}(\Gw)$ in $L^2(\Gw)$ is compact.
We proved in \cite{ChVe1}  that the positive eigenfunction   $\gg_\gm^\Gw\in H_\mu(\Gw)$ of $\CL_\gm$  associated to   the first eigenvalue $\ell_\gm^\Gw$  satisfies
\bel{1-A22}
\left.\BA {lll}
 \CL_\gm\gg_\gm^\Gw=\ell_\gm^\Gw\gg_\gm^\Gw\quad&\text{in }\, \Gw\\[1.5mm]
 \phantom{\CL_\gm}
 \gg_\gm^\Gw=0\quad&\text{on }\prt\Gw\setminus\{0\},
\EA\right.\ee
and there exist $c_1>c_2>0$ and $\tilde c>0$  such that for all $x\in\overline\Gw\setminus\{0\}$
  \bel{1-C1}\BA {lll}
(i)\qquad\qquad\qquad c_2|x|^{\ga_+-1}\gr(x)\leq \gg_\gm^\Gw(x)\leq c_1|x|^{\ga_+-1}\gr(x),\qquad\qquad\qquad\qquad\qquad\qquad\\[2mm]

(ii)\qquad\qquad\qquad  |\nabla \gg_\gm^\Gw(x)|\leq \tilde c\myfrac{\gg_\gm^\Gw(x)}{\gr(x)}.\qquad\qquad\qquad\qquad\qquad\qquad
\EA\ee
 As it is shown in \cite{ChVe1}, the function $\gg_\gm^\Gw$ plays the role of a weight function for expressing the notion of weak solutions as it is also classical in problems with boundary singularities \cite{MaVe1,MaVe2}.  Inequality $(\ref {1-A12})$ implies the existence of the Green kernel $G^\Gw_\gm$  with corresponding Green operator $\BBG^\Gw_\gm$. The Poisson kernel $K^\Gw_\gm$ of $\CL_\gm$ in $\Gw\ti\prt\Gw$ is constructed in \cite{ChVe1},  by a simple truncation as in \cite{VeYa} if $\gm\geq 0$, and by a more elaborate approximation in the general case.
When $\gm>0$ the Poisson kernel has the property that 
\bel{1-A24}
 K^\Gw_\gm(x,0)=0\quad\text{for all }\,x\in\overline\Gw\setminus\{0\},
 \ee
by \cite[Theorem A.1]{VeYa}. 
The singular kernel $\phi^\Gw_\gm$ (see \cite[Section 4]{ChVe1} for the construction) is the analogue in a bounded domain of the explicit singular solution 
$x\mapsto \phi_\gm(x)=\abs x^{\ga_--1}x_N$ defined in $\BBR^N_+$. The main property of this singular kernel is that it satisfies for all $x\in\overline\Gw\setminus\{0\}$,
  \bel{1-A25}\BA {lll}
c_3|x|^{\ga_--1}\gr(x)\leq \phi^\Gw_\gm(x)\leq c_4|x|^{\ga_--1}\gr(x) \quad \text{ if $\gm>\gm_1$},
\EA\ee
 and 
  \bel{1-A26}\BA {lll}
c_5|x|^{-\frac{N}{2}}(|\ln|x||+1)\gr(x)\leq \phi^\Gw_{\gm_1}(x)\leq c_6|x|^{-\frac{N}{2}}(|\ln|x||+1)\gr(x).
\EA\ee

We assume  in the sequel  that $\Omega$ is a bounded smooth domain such that $0\in\prt\Omega$ and that the normal vector to $\prt\Gw$ at origin is $e_N=(0,\cdots,0,1)$. We define the $\gg_\gm^\Gw$-dual operator $\CL^*_\gm$ of $\CL_\gm$ by 
 \bel{1-A27}
\CL^*_\gm\gz=-\Gd\gz-\myfrac{2}{\gg_\gm^\Gw}\langle\nabla\gg_\gm^\Gw,\nabla\gz\rangle+\ell_\gm^\Gw\gz\quad\text{for all }\,\gz\in C^{1,1}(\Gw).
 \ee
It satisfies the following commutating property
\bel{2-0B0}
 \BA {lll}
\CL_\gm(\gg^\Gw_\gm\gz)=\gg^\Gw_\gm\CL^*_\gm\gz.
\EA
\ee
We denote by $\frak M(\Gw;\gg_{\gm}^\Gw)$ the set of Radon measures $\gn$ in $\Gw$ such that
\bel{1-A28}
 \BA {lll}
 \sup\left\{\myint{\Gw}{}\gz d|\gn|:\gz\in C_c(\Gw),\,0\leq\gz\leq\gg_{\gm}^\Gw \right\}:=\myint{\Gw}{}\gg_{\gm}^\Gw d|\gn|<\infty.
\EA\ee
Thus, if $\gn\in \frak M_+(\Gw;\gg^\Gw_\gm)$ the measure $\gg^\Gw_\gm\gn$ is a bounded measure in $\Gw$. We also set
\bel{1-A29}
\gb^\Gw_\gm(x)=-\frac{\prt\gg_{\gm}^\Gw}{\prt{\bf n}}\Big\lfloor_{\prt\Gw}. 
\ee
The space of Radon measures $\gl$ on $\prt\Gw\setminus\{0\}$ such that 
\bel{1-A30}
\BA {lll}
 \sup\left\{\myint{\prt\Gw\setminus\{0\}}{}\gz d|\gl|:\gz\in C_c(\prt\Gw\setminus\{0\}),\,0\leq\gz\leq \gb^\Gw_\gm \right\}:=
\myint{\prt\Gw\setminus\{0\}}{} \gb_\gm^\Gw d|\gl|<\infty,
\EA\ee
is denoted by $\frak M(\prt\Gw;\gb^\Gw_\gm)$. The extension of $\gl\in \frak M_+(\prt\Gw;\gb^\Gw_\gm)$  as a measure $ \gb^\Gw_\gm\gl$ in $\prt\Gw$ is given by
 \bel{1-A31}
 \BA {lll}
\myint{\prt\Gw}{}\gz d( \gb^\Gw_\gm\gl)= \sup\left\{\myint{\prt\Gw}{}\gu  \gb^\Gw_\gm d\gl:\gu\in C_c(\prt\Gw\setminus\{0\}),\,0\leq\gu\leq\gz \right\}
\quad\text{for all }\,\gz\in C(\prt\Gw)\,,\;\gz\geq 0,
\EA\ee
and by $ \gb_\gm^\Gw\gl= \gb_\gm^\Gw\gl_+- \gb_\gm^\Gw\gl_-$ if $\gl$ is a signed measure in $\frak M(\prt\Gw; \gb^\Gw_\gm)$, and we use the same notation $\frak M(\prt\Gw; \gb^\Gw_\gm)$ for the set of all such extensions. The Dirac mass at $0$ does not belong  to $\frak M(\prt\Gw; \gb^\Gw_\gm)$, but it is the limit of sequences of measures in this space.  
We proved in \cite{ChVe1} that  if $\gn\in \frak M_+(\Gw;\gg_{\gm}^\Gw)$, $\gl\in \frak M(\prt\Gw; \gb^\Gw_\gm)$ and $k\in\BBR$, the function 
\bel{1-A32*}
 u=\BBG^\Gw_\gm[\gn]+\BBK^\Gw_\gm[\gl]+k\gf_\gm^\Gw:=\BBH^\Gw_\gm[(\gn,\gl+k\gd_0)],
\ee
 is the unique function belonging to $L^1(\Gw,\gr^{-1}d\gg_\gm^\Gw)$ satisfying 
 \bel{1-A32}
\myint{\Gw}{}u\CL^*_\gm\gz d\gg_{\gm}^\Gw=\myint{\Gw}{}\gz d(\gg^\Gw_\gm\gn)+\myint{\prt\Gw}{}\gz d(\gb^\Gw_\gm\gl)+kc_\mu \gz(0),
\ee
for all $\gz\in \BBX_\gm(\Gw)=\left\{\gz\in C(\overline\Gw)\text{ s.t. }\gg_\gm^\Gw\gz\in H_\gm(\Gw)\text{ and }\,\gr\CL^*_\gm\gz\in L^\infty(\Gw)\right\}$,
where   
$$c_\gm=\left\{\BA{lll}2\sqrt{\gm-\gm_1}\myint{\BBS^{N-1}_+}{}\phi_1^2dS&\quad\text{if }\;\gm>\gm_1,
\\[3mm]
\left(\frac N2-1\right)\myint{\BBS^{N-1}_+}{}\phi_1^2dS&\quad\text{if }\;\gm=\gm_1,
\EA
\right.
$$
and where $\phi_1$ is the positive eigenfunction of $\Delta_{\BBS^{N-1}}$ (normalized by $\sup\phi_1=0$) in $H^1_0(\BBS^{N-1}_+)$ where $\BBS^{N-1}_+:=\{(x',x_N)\in\R^N:\, |x|=1,\, x_N>0\}$.
\medskip

Let $g:\BBR\mapsto\BBR$ be a continuous nondecreasing function satisfying $rg(r)\geq 0$. Thanks to this result we can construct weak solutions to the problem
\bel{1-B1}\left.\BA {lll}
\CL_\gm u+g(u)=\gn\quad&\text{in }\,\Gw\\[1mm]
\phantom{\CL_\gm +g(u)}
u=\gl+k\delta_0\quad&\text{on }\prt\Gw.
\EA\right.\ee

\begin{definition}\label{weak} Let $(\gn,\gl)\in \frak M(\Gw;\gg^\Gw_\gm)\ti\frak M(\prt\Gw; \gb^\Gw_\gm)$ and $k\in\BBR$.
A function $u\in L^1(\Gw,\gr^{-1}d\gg_{\gm}^\Gw)$ is a weak solution of (\ref{1-B1}) if $g(u)\in L^1(\Gw,d\gg_{\gm}^\Gw)$ and 
 \bel{1-B3}
 \BA {lll}
\myint{\Gw}{}\left(u\CL^*_\gm\gz +g(u)\gz\right)d\gg_\gm=\myint{\Gw}{}\gz d(\gg_\gm\gn)+\myint{\prt\Gw}{}\gz d(\gb_\gm^\Gw\gl)+kc_\mu \gz(0)\ \text{ for any $\gz\in \BBX_\gm(\Gw)$.}
\EA\ee
\end{definition}
We introduce two new specific exponents,
  \bel{1-B4}\BA {lll}
p^*_\gm=1-\myfrac{2}{\ga_-}=\myfrac{N+2+2\sqrt{\gm-\gm_1}}{N-2+2\sqrt{\gm-\gm_1}}\;\text{and }\; p^{**}_\gm=1-\myfrac{2}{\ga_+}=\myfrac{N+2-2\sqrt{\gm-\gm_1}}{N-2-2\sqrt{\gm-\gm_1}}.
\EA\ee
Note that $p^{**}_\gm$ is defined only if $N\geq 3$ and $-\frac {N^2}{4}\leq\gm<1-N$. Furthermore $p^*_0=\frac{N+1}{N-1}$ and $p^*_{\mu_1}=\frac{N+2}{N-2}$. \medskip

In the present article our first result deals with the existence of a solution with an isolated singularity on boundary: \medskip

\nind{\bf Theorem A} {\it Assume $N\geq 3$ and $\gm\geq\gm_1$, or $N=2$ and  $\gm>\gm_1$,  and let $g:\BBR\mapsto\BBR$ be a continuous nondecreasing function such that $rg(r)\geq 0$. If there holds
 \bel{1-B5}\BA {lll}
\myint{1}{\infty}\left(g(s)-g(-s)\right)s^{-1-p^*_\gm}ds<\infty\ \text{ if $\gm>\gm_1$,}
\EA\ee
or 
 \bel{1-B6}\BA {lll}
\myint{1}{\infty}\left(g(s\ln s)-g(-s\ln|s|)\right)s^{-1-p^*_{\gm_1}}ds<\infty\ \text{ if $\gm=\gm_1$,}
\EA\ee
 then for any $k\in\BBR$ there exists a unique weak solution $u_{k\gd_0}$ to 
\bel{1-B7}
\left.\BA {lll}
\CL_\gm u+g(u)=0\quad&\text{in }\,\Gw\\[1mm]
\phantom{\CL_\gm+g(u)}
u=k\gd_0\quad&\text{on }\prt\Gw.
\EA\right.\ee
Furthermore,
\bel{1-B8*}
 \BA {lll}\displaystyle
\lim_{x\to 0}\myfrac{u_{k\gd_0}(x)}{\phi_\gm^\Gw(x)}=\myfrac{k}{c_\gm}.
\EA\ee
}\medskip

When the measures do not charge the point $0$, we have a result which is similar as the one proved in \cite{GmVe}.
\medskip

\nind{\bf Theorem B} {\it Assume $N\geq 3$ and $\gm\geq\gm_1$, or $N=2$ and  $\gm>\gm_1$, and let $g:\BBR\mapsto\BBR$ be a continuous nondecreasing function such that $rg(r)\geq 0$ satisfying
 \bel{1-B8-1}\BA {lll}
\myint{1}{\infty}\left(g(s)-g(-s)\right)s^{-1-p^*_0}ds<\infty.
\EA\ee
Then for any $(\gn,\gl)\in \frak M(\Gw;\gg^\Gw_\gm)\ti\frak M(\prt\Gw;\gb^\Gw_\gm)$ there exists a unique weak solution $u$ to 
\bel{1-B9}
 \left.\BA {lll}
\CL_\gm u+g(u)=\gn\ \ &\text{in }\,\Gw\\[1mm]
\phantom{\CL_\gm+g(u)}
u=\gl\ \ &\text{on }\prt\Gw.
\EA\right.\ee
}\medskip

Finally we construct a solution to (\ref{1-B1}) without restriction on the measures by gluing solutions corresponding to Theorems A and B {\it provided $g$ satisfies the weak $\Gd_2$-condition}, a condition already introduced in \cite{ChVe}: \smallskip

\nind{\it There exists a continuous nondecreasing positive function $K:\BBR_+\mapsto\BBR_+$ such that 
  \bel{1-B10}\BA {lll}
|g(s+r)|\leq K(|r|)\left(|g(s)|+|g(r)|\right)\quad\text{for all }\;(s,r)\in\BBR\ti\BBR\;\text{ s.t. }\;sr\geq 0.
\EA\ee}
\smallskip

\nind{\bf Theorem C} {\it Assume $N\geq 3$ and $\gm\geq\gm_1$, or $N=2$ and  $\gm>\gm_1$, and let $g:\BBR\mapsto\BBR$ be a continuous nondecreasing function such that $rg(r)\geq 0$ satisfying the weak $\Gd_2$-condition and 
 \bel{1-B12}\BA {lll}
\myint{1}{\infty}\left(g(s)-g(-s)\right)s^{-1-\min\{p^*_\gm,p^*_0\}}ds<+\infty.
\EA\ee
Then for any $(\gn,\gl)\in \frak M(\Gw;\gg^\Gw_\gm)\ti\frak M(\prt\Gw; \gb^\Gw_\gm)$ and $k\in\BBR$ there exists a solution $u$ to the problem (\ref{1-B1}).}\medskip

A nonlinearity $g$ for which problem $(\ref{1-B1})$ admits a solution is called {\it subcritical}. A couple of measures $(\gn,\gl)$ for which problem $(\ref{1-B1})$ admits a solution is called {\it $g$-good} (see \cite{BMP} in the case $\gm=0$). In the supercritical case all the measures are not $g$-good. Besides the problem at $0$ where 
(\ref{1-B5})-(\ref{1-B6}) may or may not be satisfied, the admissibility of a measure depends on its concentration expressed in terms of Bessel capacities (see Adams-Hedberg book \cite{AH} for a treatment of these notions). We denote these capacities by $Cap^{^{_{\BBR^d}}}_{\ga,q}$ where $d=N$ or $N-1$. In this framework we consider only the case 
 where $g(r)=g_p(r):=|r|^{p-1}r$ with $p>1$. The following theorem is proved. \medskip
 
 \nind{\bf Theorem D} {\it Assume $\gm\geq\gm_1$ and $p>1$. \\
 1- A measure $\gn\in \frak M(\Gw;\gg^\Gw_\gm)$ is $g_p$-good if and only if it is absolutely continuous with respect to the $Cap^{^{_{\BBR^N}}}_{2,p'}$-capacity.\\
  2- A measure  $\gl\in \frak M(\prt\Gw; \gb^\Gw_\gm)$ is $g_p$-good if and only if it is absolutely continuous with respect to the $Cap^{^{_{\BBR^{N-1}}}}_{\frac{2}{p},p'}$-capacity.
}\medskip

Similarly we have a characterization of removable singularities.

\medskip
 \nind{\bf Theorem E} {\it Assume $\gm\geq\gm_1$, $p>1$ and $K\subset\overline\Gw$ is compact. Then any weak solution of 
  \bel{1-B13}\left.\BA {lll}
\CL_\gm u+g_p(u)=0&\quad\text{in }\;\Gw\cap K^c\\[1mm]
\phantom{\CL_\gm +g_p(u)}
u=0&\quad\text{on }\;\prt\Gw\cap K^c,
\EA\right.\ee
 can be extended as a solution of the same equation in $\Gw$ vanishing on $\prt\Gw$ if and only if\\
 (i) $Cap^{^{_{\BBR^{N}}}}_{2,p'}(K)=0$  if $K\subset\Gw$. \\
  (ii) $Cap^{^{_{\BBR^{N-1}}}}_{\frac 2p,p'}(K)=0$  if $K\subset\prt\Gw\setminus\{0\}$. \\
   (iii) $Cap^{^{_{\BBR^N}}}_{2,p'}(K)=0$ and $Cap^{^{_{\BBR^{N-1}}}}_{\frac 2p,p'}(K\cap\prt\Gw)$ if $K\subset\overline\Gw\setminus\{0\}$. \\
  (iv) $Cap^{^{_{\BBR^{N-1}}}}_{\frac 2p,p'}(K)=0$ and $p\geq p^*_\gm$  if $0\in K\subset\prt\Gw$ and $K\setminus\{0\}\neq\{\emptyset\}$. \\
 (v) $Cap^{^{_{\BBR^N}}}_{2,p'}(K\cap\Gw)=0$, $Cap^{^{_{\BBR^{N-1}}}}_{\frac 2p,p'}(K\cap\prt\Gw)=0$ and $p\geq p^*_\gm$ if $0\in K\subset\overline\Gw$ and $K\cap\Gw\neq\{\emptyset\}$. 
 }\medskip

At end we characterize the behaviour of solutions of 
  \bel{1-B14}\left.\BA {lll}
\CL_\gm u+g_p(u)=0&\quad\text{in }\;\Gw\\[1mm]
\phantom{\CL_\gm +g_p(u)}
u=h&\quad\text{on }\;\prt\Gw\setminus\{0\},
\EA\right.\ee
where $h\in C^3(\prt\Gw)$. When $p\geq p_\gm^*$ we prove that $u$ is indeed the very weak solution of 
  \bel{1-B14h}\left.\BA {lll}
\CL_\gm u+g_p(u)=0&\quad\text{in }\;\Gw\\[1mm]
\phantom{\CL_\gm +g_p(u)}
u=h&\quad\text{on }\;\prt\Gw.
\EA\right.\ee
The techniques we use are generalization of the ones  developed for description of internal singularities in \cite{GuVe}, and for boundary singularities with $\gm=0$  in \cite{GmVe}.  For this task, we associate a problem on $\BBS^{N-1}_+$:
  \bel{1-B15}\left.\BA {lll}
-\Gd'\gw+\left(\Gl_{p,N}+\gm\right)\gw+g_p(\gw)=0&\quad\text{in }\ \BBS^{N-1}_+\\[1mm]
\phantom{-\Gd'+\left(\Gl_{p,N}+\gm\right)\gw+g_p(\gw)}
\gw=0&\quad\text{on }\;\prt \BBS^{N-1}_+,
\EA\right.\ee
where 
  \bel{1-B16}
  \BA {lll}
\Gl_{p,N}=\myfrac{2}{p-1}\left(N-\myfrac{2p}{p-1}\right).
\EA\ee
Let $\CS_{\gm,p}$ (resp. $\CS^+_{\gm,p}$) denote the set of solutions (resp. positive solutions) of (\ref{1-B15}). We set
  \bel{1-B17}\BA {lll}
\tilde p_\gm^*=1+\myfrac{2}{a_-}:=\myfrac{N+2+2\sqrt{\gm-\gm_2}}{N-2+2\sqrt{\gm-\gm_2}}\;\text{ and }\;
\tilde p_\gm^{**}=1+\myfrac{2}{a_+}:=\myfrac{N+2-2\sqrt{\gm-\gm_2}}{N-2-2\sqrt{\gm-\gm_2}},
\EA\ee
where $\gm_2=-\left(\frac{N+2}{2}\right)^2$. Note that $\tilde p_\gm^{**}$ is defined only if $N\geq 9$ and $-\frac{N^2}{4}\leq\gm<-2N$. The role of the numbers $a_+$ and $a_-$, will be explained in the proof of the theorem.  Then we have\medskip

 \nind{\bf Theorem F} {\it Assume $\gm\geq\gm_1$ and $p>1$. \smallskip
 
 \nind 1-  $\CS_{\gm,p}$ is not reduced to $\{0\}$ if and only if $\Gl_{p,N}+\gm+N-1<0$, that is\\
(i) either $1<p<p_\gm^*$,\\
(ii) or $N\geq 3$, $\gm_1\leq\gm<1-N$ and $p>p^{**}_\gm$.\smallskip
 
 \nind 2- If $\CS^+_{\gm,p}$ is non-empty, it is reduced to one element $\gw_{\gm}$.\smallskip
 
 \nind 3- All the elements of $\CS_{\gm,p}$ have constant sign  if $\Gl_{p,N}+\gm+N-1<\Gl_{p,N}+\gm+2N\leq 0$, 
 that is:\\ 
(i) when $\gm\geq 1-N$ and $\tilde p^*_\gm\leq p< p^*_\gm$,\\
(ii) when $N\geq 3$, $-2N\leq \gm< 1-N$ and either $\tilde p^*_\gm\leq p< p^*_\gm$ or
$p^{**}_\gm<p$,\\
(iii) when $N\geq 9$ and $\gm_1\leq \gm< -2N$ and either $\tilde p^*_\gm\leq p< p^*_\gm$ or $p^{**}_\gm<p\leq \tilde p^{**}_\gm$.
}\medskip

Since any solution of (\ref{1-B14}) satisfies
    \bel{1-B18}
|u(x)|\leq c_7\gr(x)|x|^{-\frac{p+1}{p-1}}\quad\text{ for all }\;x\in \Gw\cap B_{r_0},
\ee
for some $r_0>0$ and $c_7>0$ depending on $N,p$ and $\Gw$, using a diffeomorphism, we flatten the boundary as in \cite{GmVe}, define the new function $\tilde u(y)$ by this change of variable. We set $v(t,\gs)=r^{\frac{2}{p-1}}\tilde u(r,\gs)$ with $t=\ln r$ and study the limit set $\CE_v$ of the new equation satisfied by $v(t,.)$ when $t\to-\infty$. This limit set is a connected compact subset of  $\CE_\gm$. If $u\geq 0$, $\CE_v\subset \CE_\gm^+$. Thus we prove the following.\medskip

 \nind{\bf Theorem G} {\it Assume $\gm\geq\gm_1$, $h\in C^3(\prt\Gw)$ and $u\in C^2(\Gw)\cap C(\overline\Gw\setminus\{0\})$ is a nonnegative solution of (\ref{1-B14}). If $1<p<p_\gm^*$ then\smallskip

  \nind  (i) either 
   \bel{1-B19}\!\!\!\!\!\!\!\!
\lim_{\tiny\BA {cr}\Gw\ni x\to 0\\
\;\;\frac{x}{|x|}\to\gs\in S^{N-1}_+\EA}\!\!\!\!\!\!\!\!|x|^{\frac{2}{p-1}}u(x)=\gw_\gm(\gs),
\ee
 (ii) or there exists $\ell> 0$ such that 
    \bel{1-B20}
u(x)=\ell K^\Gw_\gm(x,0)(1+o(1))\quad\text{as }\, x\in\Gw,\,x\to 0,
\ee
and $u$ is the weak solution of 
  \bel{1-B14hk}\left.\BA {lll}
\CL_\gm u+g_p(u)=0&\quad\text{in }\;\Gw\\[1mm]
\phantom{\CL_\gm +g_p(u)}
u=h+c\ell\gd_0&\quad\text{on }\;\prt\Gw.
\EA\right.\ee

 }\medskip

When $u$ is a signed solution, the situation is more delicate and we obtain only partial results. \medskip

 \nind{\bf Theorem H} {\it Assume $\gm\geq\gm_1$, $h\in C^3(\prt\Gw)$ and $u\in C^2(\Gw)\cap C(\overline\Gw\setminus\{0\})$ is a solution of (\ref{1-B14}). If $\tilde p^*_\gm\leq p<p^*_\gm$, then \\
 (a) either 
   \bel{1-B21}\!\!\!\!\!\!\!\!
\lim_{\tiny\BA {cr}\Gw\ni x\to 0\\
\;\;\frac{x}{|x|}\to\gs\in \BBS^{N-1}_+\EA}\!\!\!\!\!\!\!\!|x|^{\frac{2}{p-1}}u(x)=\pm \gw_\gm(\gs),
\ee
 (b)  or
 
 \bel{1-B21'}\!\!\!\!\!\!\!\!
\lim_{\tiny\BA {cr}\Gw\ni x\to 0\\
\;\;\frac{x}{|x|}\to\gs\in \BBS^{N-1}_+\EA}\!\!\!\!\!\!\!\!|x|^{\frac{2}{p-1}}u(x)=0.
\ee
 If we assume furthermore that $\tilde p^*_\gm< p$ and (\ref{1-B21'}) is verified, then 
there exists $\ell\in\BBR$ such that (\ref{1-B20}) and (\ref{1-B14hk}) hold.
  }\medskip
  
  In two cases the limit set is reduced to a single element of $\CE_\gm$, whatever is the structure of this set.
  
  \medskip

 \nind{\bf Theorem I} {\it Assume $\gm\geq\gm_1$, $h\in C^3(\prt\Gw)$ and $u\in C^2(\Gw)\cap C(\overline\Gw\setminus\{0\})$ is a solution of (\ref{1-B14})). \smallskip
 
\nind 1- If $N+2\sqrt{\gm-\gm_1}<4$ and $p=3$, then there exists $\gw\in \CS_{\gm,p}$ such that 
\bel{1-B25}
   \!\!\!\!\!\!\!\!
\lim_{\tiny\BA {cr}\Gw\ni x\to 0\\
\;\;\frac{x}{|x|}\to\gs\in \BBS^{N-1}_+\EA}\!\!\!\!\!\!\!\!|x|^{\frac{2}{p-1}}u(x)=\gw(\gs).\\
\ee
\nind 2- If $N=2$ and $1<p<1+\frac{2}{\sqrt{\gm+1}}$, then
  \bel{1-B23}
  \!\!\!\!\!\!\!\!
\lim_{\tiny\BA {cr}\Gw\ni x\to 0\\
\;\;\frac{x}{|x|}\to\gs\in \BBS^{1}_+\EA}\!\!\!\!\!\!\!\!|x|^{\frac{2}{p-1}}u(x)=\gw(\gs),
\ee
where $\gw$ is a solution of 
   \bel{1-B24}\left.\BA {lll}\!\!\!\!\!\!
-\gw''+\Big(\gm-\Big(\myfrac{2}{p-1}\Big)^2\Big)\gw+g_p(\gw)=0\quad\text{on }\,(0,\gp)\\[3mm]
\qquad\qquad\qquad\quad\quad \gw(0)=\gw(\gp)=0.
\EA\right.\ee
Furthermore, if $\prt\Gw$ is locally a straigh line near $0$ and the limit in (\ref{1-B24}) is zero, there exists $\ell\in\BBR$ such that (\ref{1-B20}) holds.}\medskip

We end this article with a removability result. 

  \medskip

 \nind{\bf Theorem J} {\it Assume $\gm\geq \gm_1$, $p\geq p_\gm^*$, $h\in C^3(\prt\Gw)$ and $u\in C^2(\Gw)\cap C(\overline\Gw\setminus\{0\})$ is a solution of (\ref{1-B14}). Then $u$ is actually the weak solution of (\ref{1-B14h}).}\medskip

 
 The rest of this paper is organized as follows. In section 2, we recall Kato's inequality and prove the existence and uniqueness of solutions of semilinear elliptic equation with measures sources when the nonlinearity is  subcritical.  Section 3 is devoted to the supercritical case by connecting the measures to Bessel capacities. In section 4 we give an abridged proof of the behaviour of solutions near the singular point $0$.


\mysection{The subcritical case}
\subsection{Kato inequality}

\bprop{comp} Let $N\geq 2$, $\gm\geq \gm_1$ and $g:\Gw\ti\BBR\mapsto\BBR$ be a continuous function satisfying 
$$g(s_1,x)\geq g(s_2,x)\quad\text{if }\;x\in\BBR^N_+\,\;\text{and }\;s_1\geq s_2.
$$
If $u$ and $v$ belong to $C^{1,1}(\Gw)\cap C(\overline\Gw\setminus\{0\})$ satisfy
\bel{2-A1}
\left.\BA {lll}
\CL_\gm u+g(x,u)\geq \CL_\gm v+g(x,v)&\quad\text{in }\ \Gw\\[1mm]
\phantom{\CL_\gm +g(x,u)}
u\geq v&\quad\text{on }\;\prt\Gw\setminus\{0\},
\EA\right.
\ee
and 
\bel{2-A2}
\BA {lll}
\displaystyle \liminf_{r\to 0}\;\sup_{\tiny\BA{lll}x\ni\Gw\\
|x|=r\EA}\myfrac{v(x)-u(x)}{\phi^\Gw_\gm(x)}\leq 0,
\EA\ee
then $v\leq u$ in $\Gw$. 
\es
{\bf Proof.} Set $w=v-u$, then $\CL_\gm w+h(x)w=0$ where
$$h(x)=\left\{\BA {lll}
\myfrac{g(x,v)-g(x,u)}{w}&\quad\text{if }\;w\neq 0\\[1mm]
0&\quad\text{if }\;w= 0.
\EA\right.
$$ 
Hence $h\geq 0$. For $\ge>0$, we set  $W_\ge=v-u-\ge \phi^\Gw_\gm$. Then $W_\ge\in C_c^{0,1}(\overline\Gw\setminus\{0\})$. 
There exists a sequence $\{r_n\}$ tending to $0$ such that 
$$W_\ge(x)<0\quad\text{for }\;|x|=r_n,
$$
and there holds  
$$-\Gd W_\ge+\myfrac{\gm}{|x|^2}W_\ge+hW_\ge\leq 0.
$$
Multiplying by $(W_{\ge})_+:=\max\{0,W_{\ge}\}$ and integrating yields, since $(W_{\ge})_+\in C^{1,1}_c(\overline\Gw\setminus\{0\})$
$$\myint{\Gw\setminus B_{r_n}}{}\left(\abs{\nabla (W_{\ge})_+}^2+\myfrac{\gm_1}{|x|^2}(W_{\ge})_+^2\right) dx\leq 0.
$$
Hence $(W_{\ge})_+=0$ in $\Gw\setminus B_{r_n}$, we get the result by letting $r_n\to 0$ first and then $\ge\to 0$.
\qeda

\medskip

The following form of Kato's inequality  for Schr\" odinger operators with Hardy-Leray potential 
with boundary singularity singularity is important in our approach of the concept of weak solutions to (\ref{1-B1}).
\bprop{Kato} \cite[Lemma 3.1]{ChVe1}  Assume $N\geq 3$ and $\gm\geq\gm_1$, or $N=2$ and $\gm>\gm_1$. Then for any $(f,h)\in L^1(\Gw,d\gg^\Gw_\gm)\ti L^1(\prt\Gw,d\gb^\Gw_\gm)$ there exists a unique function $u\in L^1(\Gw,|x|^{-1}d\gg_\gm^\Gw)$ satisfying 
\bel{2-A3}
\BA {lll}
\myint{\Gw}{}u\CL^*_\gm\gz\, d\gg^\Gw_\gm =\myint{\Gw}{}\gz f\,d\gg^\Gw_\gm +\myint{\prt\Gw}{}h \,d\gb^\Gw_\gm \ \, \text{ for all $\gz\in \BBX_\gm(\Gw)$.}
\EA\ee
Furthermore, for any $\gz\in \BBX^+_\gm(\Gw)=\left\{\gz\in \BBX_\gm(\Gw): \gz\geq 0\right\}$, there holds
\bel{2-A5}
\BA {lll}
\myint{\Gw}{}|u|\CL^*_\gm\gz d\gg^\Gw_\gm(x)\leq \myint{\Gw}{}\gz f{\rm sgn}(u)d\gg_\gm^\Gw(x)+\myint{\prt\Gw}{}|h|\gz d\gb_\gm^\Gw(x')
\EA\ee
and
\bel{2-A6}
\BA {lll}
\myint{\Gw}{}u_+\CL^*_\gm\gz d\gg^\Gw_\gm(x)\leq \myint{\Gw}{}\gz f{\rm sgn}_+(u)d\gg^\Gw_\gm(x)+\myint{\prt\BBR^N_+}{}h_+\gz d\gb^\Gw_\gm(x').
\EA\ee
\es

Let $\gs^\Gw_\gm\in H_\gm(\Gw) $ be the unique  variational solution of 
  \bel{2-0B1*}
 \CL_\gm\gs^\Gw_\gm= \frac{\gg^\Gw_\gm}{\min\{l^\Gw_\gm, \, \gr\}}\quad\text{in }\, \Gw\quad\ {\rm and}\quad\
 \gs^\Gw_\gm=0\quad\text{on } \prt \Gw,
\ee
 then $\gs^\Gw_\gm$ belongs to $C^{2}(\overline\Gw\setminus\{0\})$ and satisfies (see \cite[Appendix]{ChVe1})
  \bel{2-0B1}
  \BA {lll}
(i)\qquad\qquad\qquad\qquad&\gg_\gm^\Gw\leq \gs^\Gw_\gm\leq c_7\gg_\gm^\Gw\quad&\text{in }\;\Gw,\qquad\qquad\qquad\qquad\qquad\qquad\\[2mm]
(ii)\qquad\qquad&\nabla \gs^\Gw_\gm(x)\sim \nabla \gg_\gm^\Gw(x)\quad&\text{as }\;x\to 0.
\EA\ee
 Furthermore $\myfrac{\prt\gs^\Gw_\gm}{\prt{\bf n}}<0$ on $\prt\Gw\setminus\{0\}$. The function
$
\eta=\myfrac{\gs^\Gw_\gm}{\gg^\Gw_\gm}
$
 which verifies 
     \bel{2-0B3}\BA {lll}
\CL^*_\gm\eta=\myfrac{1}{\min\{l^\Gw_\gm, \, \gr\}}\quad\text{in }\,\Gw,
\EA\ee
plays an important role as a test function because of the following estimates that it satisfies
 \bel{2-0B1-e}
    1\leq \eta \leq c_7 \quad \text{ and }\quad |\nabla \eta| \leq c_7\rho^{-1}\quad   \text{in }\;\Gw.
 \ee

\subsection{Proof of Theorem A}
 Assume $\Gw\subset B_1$ and let $k>0$. If $\gm>\gm_1$, we have by $(\ref{1-C1})$ and $(\ref{1-A25})$
 \bel{2-A6-1}\BA {lll}
\myint{\Gw}{}g(k\gf_\gm)d\gg_\gm\leq c_9\myint{B_R}{}g(c_8|x|^{\ga_-})|x|^{\ga_+}dx
\leq c_{10}\myint{0}{R}g(c_8r^{\ga_-})r^{\ga_++N-1}dr\\[4mm]
\phantom{\myint{\Gw}{}g(k\gf_\gm)d\gg_\gm^\Gw}
\leq c_{11}\myint{R^{1/\ga_-}}{\infty}g(s)s^{-1+\frac{\ga_++N}{\ga_-}}ds=c_{11}\myint{R^{1/\ga_-}}{\infty}g(s)s^{-1-p^*_\gm}<\infty,
\EA\ee
where $\gf_\gm(x)=|x|^{\ga_--1}x_N\geq \gf_\gm^\Gw(x)$ in $\Gw$, and $p^*_\gm$ is defined in (\ref{1-B4}). If $\gm=\gm_1$ we obtain similarly
$$\myint{\Gw}{}g(k\gf_{\gm_1})d\gg_{\gm_1}^\Gw\leq c_{11}\myint{R^{1/\ga_-}}{\infty}g(s\ln s)s^{-\frac{2N}{N-2}}ds<\infty.
$$
  For $r>0$ small enough set $\Gw_r=\Gw\setminus\overline B_r$, $\prt\Gw_r=\Gg_{1,r}\cup \Gg_{2,r}$ where 
$\Gg_{1,r}=  B^c_r\cap\prt\Gw$ and $\Gg_{2,r}= \prt B_r\cap \Gw$. We consider the problem 
 \bel{2-A7}
 \left.\BA {lll}
\CL_\gm v+g(v)=0\quad&\text{in }\,\, \Gw_r\\[1mm]
\phantom{\CL_\gm +g(v)}
v=k\gf^\Gw_\gm\quad&\text{on }\prt\Gw_r.
\EA\right.\ee
The associated functional where $G(r)=\myint{0}{r}g(s) ds$ is expressed by
$$J_\gm^r(v)=\myint{\Gw_r}{}\left(\myfrac12|\nabla v|^2+\myfrac{\gm}{2|x|^2}v^2+G(v)\right)dx,$$
and is defined over $\CH_r=\{v\in H^1(\Gw_\ge):v=k\gf^\Gw_\gm\text{ on }\prt\Gw_r\}$. Any $v\in \CH_r$ can be written as 
$v=k\gf^\Gw_\gm+w$ where $w\in H^1_0(\Gw_r)$, then $J_\gm^r(v)=J_\gm^r(k\gf^\Gw_\gm+w)=\tilde J_\gm^r(w),$ where

$$\BA {lll}\tilde J_\gm^r(w)=\myint{\Gw_r}{}\left(\myfrac12|\nabla w|^2+\myfrac{\gm}{2|x|^2}w^2+G(w+k\gf^\Gw_\gm)\right)dx%
+\myfrac{k^2}{2}\myint{\Gw_r}{}\left(|\nabla \gf^\Gw_\gm|^2+\myfrac{\gm}{|x|^2}(\gf^\Gw_\gm)^2\right)dx
\\[4mm]\phantom{\tilde J_\gm^r(\gh)}
\quad\ +\myint{\Gw_r}{}\left(\langle\nabla\gf^\Gw_\gm,\nabla w\rangle+\myfrac{\gm}{|x|^2}\gf^\Gw_\gm w\right) dx
\\[4mm]\phantom{\tilde J_\gm^r(\gh)}
=\myint{\Gw_r}{}\left(\myfrac12|\nabla w|^2+\myfrac{\gm}{2|x|^2}w^2+G(w+k\gf^\Gw_\gm)\right)dx+\myfrac{k^2}{2}\myint{\Gw_r}{}\left(|\nabla \gf^\Gw_\gm|^2+\myfrac{\gm}{|x|^2}(\gf^\Gw_\gm)^2\right)dx\\[4mm]\phantom{\tilde J_\gm^r(\gh)}
\quad\ +\myint{\Gw_r}{}\gh\CL_\gm\gf^\Gw_\gm dx+\myint{\prt\Gw_r}{}\myfrac{\prt\gf^\Gw_\gm}{\prt{\bf n}}w dS
\\[4mm]\phantom{\tilde J_\gm^r(w)}
\geq \myfrac{1}{4}\myint{\Gw_r}{}\myfrac{w^2}{|x|^2\ln^2(|x|)}dx+\myfrac{k^2}{2}\myint{\Gw_r}{}\left(|\nabla \gf^\Gw_\gm|^2+\myfrac{\gm}{|x|^2}(\gf^\Gw_\gm)^2\right)dx,
\EA$$
since $w\in H^1_0(\Gw_r)$, $(\ref{1-A12})$ holds and $G\geq 0$. Hence $\tilde J_\gm^r$ and therefore $J_\gm^r$ is coercive and since it is convex, it admits a unique minimum $u_r$, which is the unique classical solution of (\ref{2-A7}) by standard regularity and  by \rprop{comp}   such that 
  $0<u_r\leq \gf^\Gw_\gm$ in $\Gw_r$.

  By monotonicity $0<u_r\leq u_{r'}$ in $\Gw_{r'}$ if $r\in(0,\,r')$. Let $u_k=\lim_{r\to 0}u_r$. Because of $(\ref{2-A6-1})$, 
$g(u_r)\to g(u_0)$ in $L^1(\Gw;d\gg^\Gw_\gm)$.  Let $\gg_r:=\gg_\gm^{\Gw_r}$ be the first eigenfuntion of the operator 
$$\gw\mapsto-\Gd\gw+\myfrac{\gm}{|x|^2}\gw\quad \text{ in $H^1_0(\Gw_r)$}
$$
with corresponding eigenvalue $\ell_r:=\ell_\gm^{\Gw_r}$. We normalize $\gg_r$ by $\gg_r(x_0)=1$ for some fixed 
$x_0$ in $\Gw_{\frac14}$. Then $\ell_r>\ell_\gm^{\Gw}$ and $\ell_r\to\ell_\gm^{\Gw}$ when $r\to 0$. Furthermore $\gg_r\to \gg_\gm^{\Gw}$
 uniformly on $\overline\Gw_r$ for any $\gd>0$, where $\gg_\gm^{\Gw}(x_0)=1$. If $\gz\in\BBX_\gm(\Gw)$, we have
 \bel{2-A9}\BA {lll}
0=\myint{\Gw_r}{}\gz\gg_r\left(\CL_\gm u_r+g(u_r)\right) dx\\[4mm]
\phantom{0}
=\myint{\Gw_r}{}\left(\phantom{a^{a^b}}\!\!\!\!\!\!\!\!\!\left(-\gg_r\Gd\gz-2\langle\nabla\gg_r,\nabla\gz\rangle+\ell_r \gz\gg_r\right)u_r+\gz\gg_r g(u_r)\right)dx
-k\myint{\Gg_{2,r}}{}\gz\myfrac{\prt \gg_r}{\prt{\bf n}}\phi_\gm^\Gw dS.
\EA\ee
Since 
$$\BA {lll}
-\myint{\Gg_{2,r}}{}\myfrac{\prt \gg_r}{\prt{\bf n}}\phi_\gm^\Gw dS=\myint{\Gw_{\ge}}{}\phi_\gm^\Gw\Gd \gg_r dS
-\myint{\Gw_r}{}\gg_r\Gd\phi_\gm^\Gw  dS=-\ell_\ge\myint{\Gw_r}{}\gg_r\phi_\gm^\Gw dx,
\EA$$
then we obtain the following, by letting $r\to 0$, 
$$\displaystyle\lim_{r\to 0}\myint{\Gg_{2,r}}{}\myfrac{\prt \gg_r}{\prt{\bf n}}\phi_\gm^\Gw dS=\ell^\Gw_\gm\myint{\Gw}{}\gg^\Gw_\gm\phi_\gm^\Gw dx.
$$
Noting from $(\ref{2-A9})$ that
$$\displaystyle\lim_{r\to 0}
\myint{\Gw_r}{}\left(\phantom{a^{a^b}}\!\!\!\!\!\!\!\!\!\left(-\gg_r\Gd\gz-2\langle\nabla\gg_r,\nabla\gz\rangle+\ell_r \gz\gg_r\right)u_r+\gz\gg_r g(u_r)\right)dx=\myint{\Gw}{}\left(u_k\CL_\gm^*\gz+\gz g(u_k)\right)d\gg^\Gw_\gm,
$$
we infer
 \bel{2-A10}\BA {lll}
\myint{\Gw}{}\left(u_k\CL_\gm^*\gz+\gz g(u_k)\right)d\gg^\Gw_\gm
=c_{N,\gm,\Gw}k\gz(0),
\EA\ee
with 
 \bel{2-A11}
 \BA {lll}
c_{N,\gm,\Gw}=\ell^\Gw_\gm\myint{\Gw}{}\gg^\Gw_\gm\phi_\gm dx.
\EA\ee
Since $x\mapsto k\gf^\Gw_\gm(x)$ satisfies $(\ref{2-A7})$ with $g=0$, it satisfies also $(\ref{2-A10})$, always with $g=0$. Combining this result with the uniqueness and the estimates given in \cite[Proposition 2.1]{ChVe1}, we can compute the explicit value of   $c_{N,\gm,\Gw}=c_\mu$.
\qeda\smallskip

\subsection{Proof of Theorem B }

  We first assume that $(\gn,\gl)\in \frak M_+(\Gw;\gg^\Gw_\gm)\ti\frak M_+(\prt\Gw; \gb^\Gw_\gm)$.  Since $g$ satisfies $(\ref{1-A7})$ and $\CL_\mu$ is uniformly elliptic in $\Omega_r$,  it follows from \cite[Section 3]{Ve-HandBook}  that   the problem 
 \bel{2-B1-r}\left.\BA {lll}
\CL_\mu u+g(u)=\gn_\ge\quad&\text{in }\;\Gw_r\\[1mm]
\phantom{\CL_\mu  +g(u) }
u=\gl_\ge\quad&\text{on }\;\Gg_{1,r}:=\prt\Gw\cap B^c_r\\[0.8mm]
\phantom{\CL_\mu  +g(u) }
u=0\quad&\text{on }\;\Gg_{2,r}:=\Gw \cap  \prt B_r,
\EA\right.\ee
admits a unique weak solution $u_{\ge,r}$, where $\gn_\ge=\gn_\ge\chi_{B_\epsilon^c}$, $\gl_\ge=\gl_\ge\chi_{B_\epsilon^c}$
and   $0<r<\ge/2$.

 By the comparison principle, for $0<\ge'<\ge$ and $0<r'<r$ there holds
 \bel{2-B2}
(i)\quad 0<u_{\ge,r}<u_{\ge',r'}\quad\text{and }\ (ii)\quad u_{\ge,r}\leq \BBG^{\Gw_r}_{\gm}[\gn_\ge]+\BBK^{\Gw_r}_{\gm}[\gl_\ge]\leq \BBG^\Gw_{\gm}[\gn]+  \BBK^\Gw_{\gm}[\gl]\quad\text{in }\;\Gw_r,
\ee
where $\BBG^{\Gw_r}_{\gm}$ and $\BBK^{\Gw_r}_{\gm}$ denote respectively the Green and the Poisson potentials of the operator $\CL_\gm$ in $\Gw_r$.
The mappings $r\mapsto u_{\ge,r}$, $r\mapsto \BBG^{\Gw_r}_{\gm}$ and $r\mapsto \BBK^{\Gw_r}_{\gm}$ are 
decreasing. We set $\displaystyle u_\ge=\lim_{r\to 0}u_{\ge,r}$, then 
 \bel{2-B3}
 0\leq u_\ge\leq \BBG^\Gw_{\gm}[\gn_\ge]+\BBK^\Gw_{\gm}[\gl_\ge]\leq \BBG^\Gw_{\gm}[\gn]+\BBK^\Gw_{\gm}[\gl].
 \ee
If $\gz\in\BBX_\gm(\Gw)$ vanishes in some neighbourhood of $0$,  there holds for $r>0$ small enough,
 \bel{2-B4}
 \myint{\Gw_r}{}\left(u_{\ge,r}\CL^*_\gm\gz+g(u_{\ge,r})\gz\right)d\gg_\gm^\Gw=\myint{\Gw_r}{}\gz d(\gg_\gm^\Gw\gn_\ge)
+\myint{\Gg_{1,\gd}}{}\gz d(\gb_\gm^\Gw\gl_\ge).
 \ee
Letting $r\to 0$, we obtain the identity
  \bel{2-B5}\myint{\Gw}{}\left(u_{\ge}\CL^*_\gm\gz+g(u_{\ge})\gz\right)d\gg_\gm^\Gw=\myint{\Gw}{}\gz d(\gg_\gm^\Gw\gn_\ge)
+\myint{\prt\Gw}{}\gz d(\gb_\gm^\Gw\gl_\ge).
 \ee
 Because $u_\ge$ is $\CL_\gm$-harmonic in $\Gw\cap B_\ge$ and vanishes on  $\prt\Gw\cap B_\ge$, it satisfies 
 $u_\ge(x)\leq c_{12}\gg_\gm^\Gw(x)$ if $x\in \Gw\cap B_{\frac\ge 2}$ for some $c_{12}>0$ depending also on $\ge$, and 
 $(\gg_\gm^\Gw(x))^{-1}u_\ge(x)\to c_{13}\geq 0$ when $x\to 0$ by \cite[Section 3]{ChVe1}. Let $\gz\in\BBX_\gm(\Gw)$ and 
\bel{2-B5a}
 \ell_n(x)=\left\{\BA{lll}0&\quad\text{if }\;|x| <\frac1 {n}\\[2mm]
 \frac{1}{2}- \frac{1}{2}\cos\left(n\gp \left(|x|-\frac{1}{n}\right)\right)&\quad\text{if }\;\frac1 {n}\leq |x|\leq\frac 2 {n}\\[2mm]
 1&\quad\text{if }\;|x|>\frac2 {n}.
 \EA\right.
\ee
 We set $\gz_n=\ell_n\gz$. Then
 \bel{2-B6}\myint{\Gw}{}\left(u_{\ge}\CL^*_\gm\gz_n+g(u_{\ge})\gz_n\right)d\gg_\gm^\Gw=\myint{\Gw}{}\gz_n d(\gg_\gm^\Gw\gn_\ge)
+\myint{\prt\Gw}{}\gz_n d(\gb_\gm^\Gw\gl_\ge).
 \ee
Firstly we observe that
$$\myint{\Gw}{}\gz_n d(\gg_\gm^\Gw\gn_\ge)
+\myint{\prt\Gw}{}\gz_n d(\gb_\gm^\Gw\gl_\ge)\to \myint{\Gw}{}\gz d(\gg_\gm^\Gw\gn_\ge)
+\myint{\prt\Gw}{}\gz d(\gb_\gm^\Gw\gl_\ge)\quad\text{as }\;n\to\infty.
$$
Then, for $n$ large enough, 
$$\myint{\Gw}{}g(u_{\ge})\gz_n d\gg_\gm^\Gw=\myint{\Gw_{\frac r2}}{}g(u_{\ge})\gz_n d\gg_\gm^\Gw+
\myint{\Gw\cap B_{\frac r2}}{}g(u_{\ge})\gz_n d\gg_\gm^\Gw=:A_n+B_n.
$$
Because $G^\Gw_{\gm}$ and $K^\Gw_{\gm}$  are respectively equivalent to $G^\Gw_{0}$ and $K^\Gw_{0}$ in $\Gw_{\frac r2}$, the condition 
(\ref{1-B8-1}), jointly with (\ref{2-B3}), implies that $A_n$ is bounded independently of $n$ and converges to $\myint{\Gw_{\frac r 2}}{}g(u_{\ge})\gz d\gg_\gm^\Gw$. If $\gm\geq 1-N$, $\ga_+$ is nonnegative thus $g(u_{\ge})\gz_n \gg_\gm^\Gw$ is bounded in $B_{\frac r2}$. If $\gm_1\leq\gm<1-N$, then $\ga_+<0$ and we have 
$$\BA {lll}|B_n|\leq \myint{\Gw\cap B_{\frac r2}}{}g(c_{12})\gz_n d\gg_\gm^\Gw\leq 
\myint{0}{r}g(c_{12}r^{\ga_+})r^{\ga_++N-1} dr\leq \myfrac{1}{\ga_+}\myint{r^{\frac{1}{\ga_+}}}{\infty}g(c_{12}s)s^{\frac{N}{\ga_+}}ds<\infty
\EA$$
 since $\frac{N}{\ga_+}\leq -1-\frac{N+2}{N-2}<-1-\frac{N+1}{N-1}$ and (\ref{1-B8-1}) holds. Therefore,
 $$\displaystyle\lim_{n\to\infty}\myint{\Gw}{}g(u_{\ge})\gz_n d\gg_\gm^\Gw=\myint{\Gw}{}g(u_{\ge})\gz d\gg_\gm^\Gw.
 $$
 Finally, we estimate the term
 $$\myint{\Gw}{}u_{\ge}\CL^*_\gm\gz_n d\gg_\gm^\Gw=C_n+D_n+E_n,
 $$
 with
 $$C_n=\myint{\Gw}{}\ell_nu_{\ge}\CL^*_\gm\gz d\gg_\gm^\Gw\;,\;\, D_n=\myint{\Gw}{}\gz u_{\ge}\CL^*_\gm\ell_n d\gg_\gm^\Gw
 \;,\;\, E_n=-2\myint{\Gw}{}u_{\ge}\langle\nabla\gz,\nabla\ell_n\rangle d\gg_\gm^\Gw. 
 $$
 Since $u_\ge$ satisfies (\ref{2-B3}) it follows from \cite[Theorem D]{ChVe1} that it is bounded in $L^1(\Gw,\gr^{-1}d\gg_\gm^\Gw)$ independently of 
 $\ge$. Hence 
 $$\displaystyle\lim_{n\to\infty} C_n=\myint{\Gw}{}u_{\ge}\CL^*_\gm\gz d\gg_\gm^\Gw.
 $$
 Using the fact that $u_\ge(x)\sim c_{13}\gg_\gm^\Gw(x)$ and $\gz(x)=\gz(0)(1+o(1))$ when $x\to 0$ we obtain
 $$\BA {lll}\myint{\Gw}{}\gz u_{\ge}\CL^*_\gm\ell_n d\gg_\gm^\Gw
 =c_{13}\gz(0)\myint{\Gw\cap \left(B_{\frac2n}\setminus B_{\frac1{n}}\right)}{}\left(-(\gg_\gm^\Gw)^2\Gd\ell_n -2\gg_\gm^\Gw\langle\nabla\ell_n,\nabla \gg_\gm^\Gw\rangle
\right)dx+o(1)=o(1),
 \EA$$
 since $\ell_n'(\frac1n)=\ell'_n(\frac2{n})=0$ and $\gg_\gm^\Gw$ vanishes on $\prt\Gw$. Similarly 
 $$\displaystyle\lim_{n\to\infty} E_n=0.
 $$
 These facts imply that 
  \bel{2-B8}\myint{\Gw}{}\left(u_{\ge}\CL^*_\gm\gz+g(u_{\ge})\gz\right)d\gg_\gm^\Gw=\myint{\Gw}{}\gz d(\gg_\gm^\Gw\gn_\ge)
+\myint{\prt\Gw}{}\gz d(\gb_\gm^\Gw\gl_\ge)\ \ \text{for any }\gz\in\BBX_\mu(\Omega).
 \ee
Notice  that from the above derivation,  (\ref{2-B8}) holds true for $\gz=\eta$, where $\eta$
is defined in (\ref{2-B8}).

 Hence $u_\ge$ is the weak solution of 
 \bel{2-B1}\left.\BA {lll}
\CL_\mu u+g(u)=\gn_\ge\quad&\text{in }\ \Gw\\[1mm]
\phantom{\CL_\mu  +g(u) }
u=\gl_\ge\quad&\text{on }  \prt\Gw.
\EA\right.\ee 
Because of uniqueness, $\ge\mapsto u_\ge$ is increasing and  $u:=\displaystyle \lim_{\ge\to 0}u_\ge$ satisfies
\bel{2-B9*}
0\leq u\leq \BBG^\Gw_{\gm}[\gn]+  \BBK^\Gw_{\gm}[\gl]\quad\text{in }\;\Gw.
\ee
If we take $\gz=\eta$ defined by (\ref{2-0B3}) we deduce from (\ref{2-B8})
\bel{2-B9}
\myint{\Gw}{}\left(\myfrac{u_\ge}{\gr}+g(u_{\ge})\eta\right)d\gg_\gm^\Gw=\myint{\Gw}{}\eta d(\gg_\gm^\Gw\gn_\ge)
+\myint{\prt\Gw}{}\eta d(\gb_\gm^\Gw\gl_\ge).
\ee
 The right-hand side of the above identity converges to $\myint{\Gw}{}\eta d(\gg_\gm^\Gw\gn)
+\myint{\prt\Gw}{}\eta d(\gb_\gm^\Gw\gl)$. Then by monotone convergence
$$
\myint{\Gw}{}\left(\myfrac{u}{\gr}+g(u)\eta\right)d\gg_\gm^\Gw=\myint{\Gw}{}\eta d(\gg_\gm^\Gw\gn)
+\myint{\prt\Gw}{}\eta d(\gb_\gm^\Gw\gl).
$$
This implies that $u_\ge\to u$ in $L^1(\Gw,\gr^{-1}d\gg_\gm^\Gw)$ and $g(u_\ge)\to g(u)$ in $L^1(\Gw,d\gg_\gm^\Gw)$ as $\ge\to0^+$. Therefore, 
since any $\gz\in \BBX_\gm(\Gw)$, satisfies $|\gz|\leq c\eta$ for some $c>0$, we infer
\bel{2-B10}
  \myint{\Gw}{}\left(u\CL^*_\gm\gz+g(u)\gz\right)d\gg_\gm^\Gw=\myint{\Gw}{}\gz d(\gg_\gm^\Gw\gn)
+\myint{\prt\Gw}{}\gz d(\gb_\gm^\Gw\gl),
\ee
 which completes the proof when the two measures are nonnegative. 
 
 In the general case we use the Jordan decomposition $\gn=\gn^+-\gn^-$, $\gl=\gl^+-\gl^-$ where $\gn^+$, $\gn^-$, $\gl^+$ and $\gl^-$ are nonnegative. Let $\gn_\ge^{\pm}$ and $\gl_\ge^{\pm}$ be $\gn^{\pm}\chi_{_{\Gw_\ge}}$ and $\gl^{\pm}\chi_{_{\prt\Gw\cap B^c_\ge}}$ respectively. We denote by $u^+_{\ge,r}$  the solution of 
 (\ref{2-B1-r}) corresponding to the couple $(\gn_\ge^{+},\gl_\ge^{+})$ and  by $u^-_{\ge,r}$ the solution of 
 \bel{2-B11}
 \left.\BA {lll}
\CL u-g(u)=\gn^-_\ge\quad&\text{in }\;\Gw_r\\[1mm]
\phantom{\CL -g(u)}
u=\gl^-_\ge\quad&\text{on }\;\Gg_{1,r} \\[1mm]
\phantom{\CL -g(u)}
u=0\quad&\text{on }\;\Gg_{2,r}.
\EA\right.\ee
Then $-u^-_{\ge,r}\leq \min\{0,u_{\ge,r}\}\leq \max\{0,u_{\ge,r}\}\leq u^+_{\ge,r}$. The mapping $r\mapsto u^+_{\ge,r}$ (resp.
$r\mapsto u^-_{\ge,r}$) is monotone increasing and we set $\displaystyle u^+_\ge=\lim_{r\to 0} u^+_{\ge,r}$ (resp. $\displaystyle u^-_\ge=\lim_{r\to 0} u^-_{\ge,r}$). The mapping $r\mapsto u_{\ge,r}$ has no reason to be monotone, but by standard regularity theory there exists $\{r_j\}$ converging to $0$ and $u_\ge\in L^{q}_{loc}$ ($1<q<\frac{N}{N-1}$) such that $u_{\ge,r_j}\to u_{\ge}$ in $L^{q}_{loc}(\Gw)$ and a.e. in $\Gw$. Hence $u_{\ge}$ satisfies (\ref{2-B5}). Since (\ref{2-B6}) holds we derive that $u_\ge$ satisfies (\ref{2-B1}). We end the proof as in the first case, using dominated convergence theorem.\qeda\smallskip

\subsection{Proof of Theorem C}
 We first assume that $\nu,\gl$ and $k$ are nonnegative. For $0<r<\ge/4$ we consider the problem
\bel{2-C1} \left. \BA {lll}
\CL_\mu u+g(u)=\gn_\ge&\quad\text{in }\ \Gw_r\\[1mm]
\phantom{\CL_\mu  +g(u)}u=\gl_\ge&\quad\text{on }\;\Gg_{1,r}\\[1mm]
\phantom{\CL_\mu  +g(u)}u=k\phi^\Gw_\mu&\quad\text{on }\;\Gg_{2,r}.
\EA\right.\ee
  The solution is denoted by $u_{\ge,k,r}$ and we recall that $u_{\ge,r}$ is the solution of 
(\ref{2-B1-r}). There holds
 \bel{2-C2}\BA {lll}
\max\{u_{\ge,r},u_{k\gd_0}\}\leq u_{\ge,k,r}\leq u_{\ge}+k\BBK^\Gw_\gm[\gd_0]\quad \text{in $\Gw_r$.}
\EA\ee
 Furthermore $u_{\ge,k,r}\leq u_{\ge,k,r'}$ if $0<r'<r$. Since $u_{\ge}$ and $k\BBK^\Gw_\gm[\gd_0]$ belong to 
$L^1(\Gw,\gr^{-1}d\gg_\gm^\Gw)$ it implies that $u_{\ge,k,r}$ converges in $L^1(\Gw,\gr^{-1}d\gg_\gm^\Gw)$ and almost everywhere 
 to $u_{\ge,k}$ when $r\to 0$. Since $\gg_\gm^\Gw$ is a supersolution for the equation $\CL_\gm u+g(u)=0$ in $\Gw_r$, for any $0<\ge_0<\ge/4$ there exists $c_{14}>0$ depending on $\ge_0$ such that for $0<r\leq\ge_0/4$, 
 $$u_{\ge,r}(x)\leq c_{14}\gg_\gm^\Gw(x)\qquad\text{for all }\; x\in B_{\ge_0}\cap\Gw_r.$$
 For any $\gs>0$ there exists $r_\gs>0$ such that for any $r<r_\gs$, $u_{\ge,r}\leq \gs \BBK^\Gw_\gm[\gd_0]$ in $B_{r_\gs}\cap\Gw_r$. Therefore 
 $u_{\ge}+k\BBK^\Gw_\gm[\gd_0]\leq (k+\gs)\BBK^\Gw_\gm[\gd_0]$ in $B_{r_\gs}\cap\Gw$. This implies
  \bel{2-C3}\BA {lll}
g(u_{\ge,k,r})\leq g((k+\gs)\BBK^\Gw_\gm[\gd_0])\qquad\text{in }\; \Gw_r\cap B_{r_\gs}.
\EA\ee
Then we obtain, with $R=\,$diam$\,\Gw$ and for some $c>0$,
$$\BA{lll}
\myint{\Gw}{}g((k+\gs)\BBK^\Gw_\gm[\gd_0])d\gg_\gm^\Gw\leq \myint{0}{R}g(c|x|^{\ga_-})|x|^{\ga_+}dx\\[4mm]
\phantom{\myint{\Gw}{}g((k+\gs)\BBK^\Gw_\gm[\gd_0])d\gg_\gm^\Gw}
=\myfrac{1}{|\ga_-|}\myint{R^{\frac{1}{\ga_-}}}{\infty}g(ct)t^{\frac{N+\ga_+-\ga_-}{\ga_-}}dt\leq \myfrac{1}{|\ga_-|}\myint{R^{\frac{1}{\ga_-}}}{\infty}g(ct)t^{-1-p^*_\gm}<\infty.
\EA$$
This implies in particular that 
  \bel{2-C4}\BA {lll}
\myint{\Gw_r\cap B_{r_\gs}}{}g((k+\gs)\BBK^\Gw_\gm[\gd_0])d\gg_\gm^\Gw\leq  \myfrac{1}{|\ga_-|}\myint{R^{\frac{1}{\ga_-}}}{\infty}g(ct)t^{-1-p^*_\gm}.
\EA\ee
In the set $\Gw_{r_\gs}$, we have $k\BBK^\Gw_\gm[\gd_0]\leq cr^{\ga_-}_\gs$ for some $c>0$. By the local $\Gd_2$-condition, we deduce
  \bel{2-C5}g(u_{\ge,r,k})\leq g(u_{\ge}+k\BBK^\Gw_\gm[\gd_0])\leq K(cr^{\ga_-}_\gs)\left(g(u_{\ge})+g(cr^{\ga_-}_\gs)\right).
\ee
Because $g(u_{\ge})$ is bounded in $L^1(\Gw_r,d\gg_\gm^{\Gw})$ independently of $r$ by Theorem B, we infer from (\ref{2-C3}), (\ref{2-C4}) and (\ref{2-C5})  that $g(u_{\ge,k,r})$ is bounded in $L^1(\Gw_r,d\gg_\gm^{\Gw})$ independently of $r$. Let $\gz\in \BBX_\gm(\Gw)$ vanishing near $0$, then for $r$ small enough,
 \bel{2-C6*}
  \myint{\Gw}{}\left(u_{\ge,k,r}\CL^*_\gm\gz+g(u_{\ge,k,r})\gz\right)d\gg_\gm^\Gw=\myint{\Gw}{}\gz d(\gg_\gm^\Gw\gn_\ge)
+\myint{\prt\Gw}{}\gz d(\gb_\gm^\Gw\gl_\ge).
\ee
Using the mononoticity of $r\mapsto u_{\ge,k,r}$  and the dominated convergence theorem we get
\bel{2-C6}
  \myint{\Gw}{}\left(u_{\ge,k}\CL^*_\gm\gz+g(u_{\ge,k})\gz\right)d\gg_\gm^\Gw=\myint{\Gw}{}\gz d(\gg_\gm^\Gw\gn_\ge)
+\myint{\prt\Gw}{}\gz d(\gb_\gm^\Gw\gl_\ge).
\ee
As we can notice it, the singular measure $k\gd_0$ cannot appear in this formulation. If $\gz\in\BBX_\gm(\Gw)$ we set $\gz_n=\ell_n\gz$ where $\ell_n$ is defined in (\ref{2-B5a}). Then
\bel{2-C7} 
  \myint{\Gw}{}\left(u_{\ge,k}\CL^*_\gm\gz+g(u_{\ge,k})\gz\right)\ell_nd\gg_\gm^\Gw- \myint{\Gw\cap\left(B_{\frac{2}{n}}\setminus B_{\frac{1}{n}}\right)}{}A_nu_{\ge,k}d\gg_\gm^\Gw=
\myint{\Gw}{}\gz \ell_nd(\gg_\gm^\Gw\gn_\ge)
+\myint{\prt\Gw}{}\gz \ell_nd(\gb_\gm^\Gw\gl_\ge),
\ee
where
$$A_n=\gz\Gd\ell_n+2\langle\nabla\ell_n,\nabla\gz\rangle+2\ga_+\gz\langle\nabla\ell_n,\frac{x}{|x|^2}\rangle.
$$
Clearly we have that 
$$\displaystyle\lim_{n\to\infty}\myint{\Gw}{}\left(u_{\ge,k}\CL^*_\gm\gz+g(u_{\ge,k})\gz\right)\ell_nd\gg_\gm^\Gw
=\myint{\Gw}{}\left(u_{\ge,k}\CL^*_\gm\gz+g(u_{\ge,k})\gz\right)d\gg_\gm^\Gw,
$$
and 
$$\displaystyle\lim_{n\to\infty}\left(\myint{\Gw}{}\gz \ell_nd(\gg_\gm^\Gw\gn_\ge)
+\myint{\prt\Gw}{}\gz \ell_nd(\gb_\gm^\Gw\gl_\ge)\right)=\myint{\Gw}{}\gz d(\gg_\gm^\Gw\gn_\ge)
+\myint{\prt\Gw}{}\gz d(\gb_\gm^\Gw\gl_\ge).
$$
Next
$$A_n\!=\!\left[\myfrac{n^2\gp^2}{2}\cos\left(\!n\gp\left(|x|-\myfrac 1n\right)\right)\!+\!\myfrac{n\gp(N-1+2\ga_+)}{2|x|}\sin\left(\!n\gp\left(|x|-\myfrac 1n\right)\right)\right]\left(\gz(0)+o(1)\right)+O(n).
$$
Using (\ref{2-C2}) with $\gd=0$ and the fact that $u_\ge=o(\BBK_\gm^\Gw[\gd_0])$ near $0$, we obtain after a technical but straightforward computation
 \bel{2-C8}  
\displaystyle\lim_{n\to\infty}\myint{\Gw\cap\left(B_{\frac{2}{n}}\setminus B_{\frac{1}{n}}\right)}{}A_nu_{\ge,k}d\gg_\gm^\Gw=kc_\gm\gz(0).
 \ee

By the normalization chosen it follows that $u_{\ge,k}$ satisfies 
 \bel{2-C9}  
\myint{\Gw}{}\left(u_{\ge,k}\CL^*_\gm\gz+g(u_{\ge,k})\gz\right)d\gg_\gm^\Gw=\myint{\Gw}{}\gz d(\gg_\gm^\Gw\gn_\ge)
+\myint{\prt\Gw}{}\gz d(\gb_\gm^\Gw\gl_\ge)+kc_\mu \gz(0).
 \ee
 Hence $u_{\ge,k}$ is the weak solution of 
 \bel{2-C10}
  \left.\BA {lll}
\CL_\mu u+g(u)=\gn_\ge&\quad\text{in }\;\Gw\\[1.5mm]
\phantom{\CL_\mu  +g(u)}u=\gl_\ge+k\gd_0&\quad\text{on }\;\prt\Gw.
\EA\right.\ee
The end of the proof in the nonnegative case is standard: we observe that the mapping $\ge\mapsto u_{\ge,k}$ is nondecreasing. We denote by $u_k$ its limit when $\ge\to 0$. If $\gz\in\BBX_\gm(\Gw)$, the right-hand side of (\ref{2-C9}) converges to 
$$\myint{\Gw}{}\gz d(\gg_\gm^\Gw\gn)
+\myint{\prt\Gw}{}\gz d(\gb_\gm^\Gw\gl)+kc_\mu \gz(0)\quad \text{as $\ge\to 0$.}
$$
 If we take $\gz=\eta$,  by property (\ref{2-0B1-e}),   (\ref{2-C8}) becomes 
 \bel{2-C8-e}  
\displaystyle\limsup_{n\to\infty}\myint{\Gw\cap\left(B_{\frac{2}{n}}\setminus B_{\frac{1}{n}}\right)}{}A_nu_{\ge,k}d\gg_\gm^\Gw\leq kc_\gm\sup_{\Omega}\eta,
 \ee
and when $\ge\to 0$,
 \bel{2-C11}  
\myint{\Gw}{}\left(\myfrac{u_k}{\gr}+g(u_{k})\eta\right)d\gg_\gm^\Gw\leq \myint{\Gw}{}\eta d(\gg_\gm^\Gw\gn_\ge)
+\myint{\prt\Gw}{}\eta d(\gb_\gm^\Gw\gl_\ge)+kc_\gm\sup_{\Omega}\eta.
 \ee
  Thus,  by the monotone convergence theorem we have that 
 $u_{\ge,k}\to u_k$ in $L^1(\Gw,\gr^{-1})d\gg_\gm^\Gw)$ and $g(u_{\ge,k})\to g(u_k)$ in $L^1(\Gw,d\gg_\gm^\Gw)$ as $\ge\to 0$. Therefore, by the dominated convergence theorem we conclude that for any $\gz\in\BBX_\gm(\Gw)$ there holds
  \bel{2-C12}  
\myint{\Gw}{}\left(u_{k}\CL^*_\gm\gz+g(u_{k})\gz\right)d\gg_\gm^\Gw=\myint{\Gw}{}\gz d(\gg_\gm^\Gw\gn)
+\myint{\prt\Gw}{}\gz d(\gb_\gm^\Gw\gl)+kc_\mu\gz(0).
\ee
 Hence $u_k$ is the weak solution of (\ref{1-B1}). 
When $\gn$ and $\gl$ are signed measures and $k$ is a real number, we use the Jordan decomposition of 
$\gn=\gn^+-\gn^-$ and $\gl=\gl^+-\gl^-$ and assume for example that $k$ is  nonnegative and we construct the solutions 
$u^+_{\ge,k,r}$ of 
 \bel{2-C14}
 \left.\BA {lll}
\CL_\mu u+g(u)=\gn^+_\ge&\quad\text{in }\;\Gw_r\\[1mm]
\phantom{\CL_\mu  +g(u)}u=\gl^+_\ge&\quad\text{on }\; \Gg_{1,r}\\[1mm]
\phantom{\CL_\mu  +g(u)}u=u_{k\gd_0} &\quad\text{on }\;\Gg_{2,r},
\EA\right.\ee
and $u^-_{\ge,\gd}$ of  
\bel{2-C15}
\left.\BA {lll}
\CL_\mu u-g(u)=\gn^-_\ge&\quad\text{in }\;\Gw_r\\[1mm]
\phantom{\CL_\mu+g(u)}u=\gl^-_\ge&\quad\text{on }\; \Gg_{1,r}\\[1mm]
\phantom{\CL_\mu+g(u)}u=0&\quad\text{on }\; \Gg_{2,r}.
\EA\right.\ee
Then the function $u_{\ge,k,r}$ of (\ref{2-C1}) satisfies  
$-u^-_{\ge,r}\leq\min\{0,u_{\ge,k,r}\}\leq\max\{0,u_{\ge,k,r}\}\leq u^+_{\ge,k,r}$. Since $u^+_{\ge,k,r}$
 is monotone with respect to $r$ with limit $u^+_{\ge,k}$, we obtain, as in the proof of {\it Theorem B}, the existence of 
 a limit $u_{\ge,k}$ of a sequence $u_{\ge,k,r_j}$, a.e. and in $L^{q}_{loc}(\Gw)$, and $u_{\ge,k}$ satisfies (\ref{2-C4}) 
 for any $\gz\in\BBX_\gm(\Gw)$ which vanishes near $0$.
 
Since  $u^-_{\ge, r}$ is $\CL_\mu$-harmonic in $\Omega\cap B_\epsilon$, $u^-_{\ge, r}=0$ on $(\prt\Omega \cap B_\epsilon)\setminus\{0\}$ and charges no Dirac mass at origin in the weak sense, then 
$$-u^-_{\ge,r}\geq -c_{15}\gg^\Gw_\mu\quad{\rm on}\quad \Gw\cap \prt B_{\frac\ge2},$$
for some $c_{15}>0$ dependent of $\epsilon$. Thus, there exists $c_{16}>0$ such that 
 $$u_{\ge,k,r}\geq u_{k\gd_0}-c_{16}\gg^\Gw_\mu:=w\qquad\text{for all }\; x\in \Gw\cap    B_{\frac\ge2}.
 $$
   Combining these estimates with (\ref{2-C2}) (applied to $u^+_{\ge,r,k}$) we obtain 
 \bel{2-C17}
 \BA {lll}
u_{k\gd_0}- c_{16}\gg_\gm^\Gw\leq u_{\ge,k,r}\leq u^+_{\ge,r}+k\BBK_\gm^\Gw[\gd_0]\leq u^+_{\ge}+k\BBK_\gm^\Gw[\gd_0]\quad\text{in }\;  \Gw\cap  B_{\frac\ge2},
\EA \ee
where $u^+_{\ge,r}$ and $u^+_{\ge}$ are the solutions of (\ref{2-B1}) with $r>0$ and $r=0$ respectively
with $\gn_\ge$ and $\gl_\ge$ replaced by $\gn^+_\ge$ and $\gl^+_\ge$. Thanks to estimate (\ref{2-C17}) we infer as in the case where $\gn_\ge$ and $\gl_\ge$ are nonnegative that $u_{\ge,k}$ satisfies  (\ref{2-C9}). We also have 
$$-u^-_{\ge}\leq\min\{0,u_{\ge,k}\}\leq\max\{0,u_{\ge,k}\}\leq u^+_{\ge,k},$$
and
$$g(-u^-_{\ge})\leq \min\{0,g(u_{\ge,k})\}\leq\max\{0,g(u_{\ge,k})\}\leq g(u^+_{\ge,k}).
$$
Then there exist a function $u_k\in L^{q}_{loc}(\Gw)$ ($1<q<\frac N{N-1}$) and a sequence $\{\ge_j\}$ converging to 
$0$ such that  $u_{\ge_j,k}\to u_k$ $L^{q}_{loc}(\Gw)$ and a.e. in $\Gw$. Since $g(u^+_{\ge,k})$ and $g(-u^-_{\ge})$ converge in $L^1(\Gw,d\gg_\gm^\Gw)$ and $u^+_{\ge,k}$ and $u^-_{\ge}$ in $L^1(\Gw,\gr^{-1}d\gg_\gm^\Gw)$, it follows that $g(u_{\ge,k})$ and $u_{\ge,k}$ endow the same properties. This is sufficient to see that (\ref{2-C9}) implies 
(\ref{2-C11}), which ends the proof.\qeda

\mysection{The supercritical case}
\subsection{Reduced measures}
We present here the notion of reduced measure which has been introduced by Brezis, Marcus and Ponce \cite{BMP} when $\gm=0$.
This notion turned out to be a very useful tool for analyzing supercritical problems. Since part of the results are simple adaptations to our framework of similar ones used in \cite{ChVe}, we will state them without detailled proofs, although the transcription to the types of measures used therein needs some precaution. We assume here that $g$ is a continuous nondecreasing function vanishing at $0$. For $\ell>0$ we set 
\bel{3-D1}
g_\ell(r)=\left\{\BA {lll}\min\{g(r),\,g(\ell)\}&\quad\text{if }\; r\geq 0\\[1mm]
\max\{g(-\ell),\, g(r)\}&\quad\text{if }\; r< 0.
\EA\right.
\ee
If $k\in\BBR_+$, and $(\gn,\gl)\in \frak M_+(\Gw;\gg^\Gw_\gm)\ti\frak M_+(\prt\Gw; \gb^\Gw_\gm)$ we denote by 
$u_{\ell}$ the solution of 
\bel{3-DX2}
 \left.\BA {lll}
\CL_\mu u+g_\ell(u)=\gn&\quad\text{in }\;\Gw\\[1mm]
\phantom{\CL_\mu+g_\ell(u)}u=\gl+k\gd_0&\quad\text{on }\;\prt\Gw.
\EA\right.\ee
Existence of $u_\ell$ comes from {\it Theorem C.} 
\bprop{P1} Let $k\in\BBR_+$, and $(\gn,\gl)\in \frak M_+(\Gw;\gg^\Gw_\gm)\ti\frak M_+(\prt\Gw; \gb^\Gw_\gm)$, then $\ell\mapsto u_\ell$ is monotone decreasing and converges to some function $u^*$ when 
$\ell\to\infty$ and there exists a real number $k^*\in [0,k]$ and two measures $(\gn^*,\gl^*)\in \frak M_+(\Gw;\gg^\Gw_\gm)\ti\frak M_+(\prt\Gw; \gb^\Gw_\gm)$ satisfying $0\leq\gn^*\leq\gn$ and $0\leq\gl^*\leq\gl$ such that $u^*$ is a weak solution of  
\bel{3-DX3}
\left.\BA {lll}
\CL_\mu u+g_\ell(u)=\gn^*&\quad\text{in }\;\Gw\\[1mm]
\phantom{\CL_\mu+g_\ell(u)}u=\gl^*+k^*\gd_0&\quad\text{on }\;\prt\Gw.
\EA\right.\ee
Furthermore the correspondence $(\gn,\gl,k)\mapsto (\gn^*,\gl^*,k^*)$ is nondecreasing.

\es
{\bf Proof.} The monotonicity is clear. By Fatou's lemma $\displaystyle u^*:=\lim_{\ell\to\infty}u_\ell$ satisfies 
$$\myint{\Gw}{}\left(u^*\CL_\gm^*\gz+g(u^*)\gz\right)d\gg_\gm^\Gw\leq  \myint{\Gw}{}\gz d(\gg_\gm^\Gw\gn)
+\myint{\prt\Gw}{}\gz d(\gb_\gm^\Gw\gl)+kc_\mu\gz(0)\quad\text{for all }\;\gz\in\BBX_\gm(\Gw),\,\gz\geq 0.
$$
The function $u^*$ is the largest subsolution of problem (\ref{1-B1}). Since the mapping 
$$\gz\mapsto \myint{\Gw}{}\left(u^*\CL_\gm^*\gz+g(u^*)\gz\right)d\gg_\gm^\Gw\quad\text{for all }\,\gz\in C^\infty_c(\Gw)$$is a positive distribution, it is a positive measure denoted by $\gn^*$. It is smaller than 
$\gn$, hence it belongs to $\frak M_+(\Gw;\gg^\Gw_\gm)$. Similarly the function $u^*$ admits a boundary trace $\gl^*$
on $\prt\Gw\setminus\{0\}$ which is a positive Radon measure smaller than $\gl$. Hence $\gl^*\in \frak M_+(\prt\Gw^*; \gb^\Gw_\gm)$. By using (\ref{1-A3}), it is extended as a measure on $\prt\Gw$, still denoted by  $\gl^*$. If 
$\gz\in\BBX_\gm(\Gw)$ vanishes near $0$, there holds
\bel{3-D4*}
\BA {lll}
\myint{\Gw}{}\left(u^*\CL_\gm^*\gz+g(u^*)\gz\right)d\gg_\gm^\Gw=\myint{\Gw}{}\gz d(\gg_\gm^\Gw\gn^*)
+\myint{\prt\Gw}{}\gz d(\gb_\gm^\Gw\gl^*).
\EA\ee
Let $v$ be the solution of 
\bel{3-D4}
 \left.\BA {lll}
\CL_\mu v+g(v)=\gn^*&\quad\text{in }\;\Gw\\[1mm]
\phantom{\CL_\mu+g(v)}v=\gl^*&\quad\text{on }\prt\Gw.
\EA\right.\ee
Existence is standard since $u^*$ exists. Furthermore $v$ is a subsolution of problem (\ref{1-B1}) hence it is smaller than $u^*$. Therefore $w=u^*-v$ is nonnegative and it satisfies
\bel{3-D5}
 \left.\BA {lll}
\CL_\mu w+g(u^*)-g(v)=0&\quad\text{in }\;\, \Gw\\[1mm]
\phantom{\CL_\mu+g(u^*)-g(v)}w=0&\quad\text{on }\;\prt\Gw\setminus\{0\}.
\EA\right.\ee
Let $\psi\in H_\mu$ be the solution of 
 \bel{3-D6}\left.\BA {lll}
\CL_\mu \psi=g(u^*)-g(v)&\quad\text{in }\;\Gw\\[1mm]
\phantom{\CL_\mu }\psi=0&\quad\text{on } \prt\Gw,
\EA\right.\ee
then $w+\psi$ is a nonnegative $\CL_\gm$-harmonic function vanishing on $\prt\Gw\setminus\{0\}$. By 
\cite[Theorem A]{ChVe1} there exists $k^*\geq 0$ such that 
$$\displaystyle \lim_{x\to 0}\myfrac{(w+\psi)(x)}{\gg_\gm^\Gw(x)}=k^*,
$$
and
$$\myint{\Gw}{}(w+\psi)\CL^*_\gm\gz d\gg_\gm^\Gw=k^*c_\mu\gz(0)\qquad\text{for all }\,\gz\in\BBX_\gm(\Gw).
$$
It follows from (\ref{3-D6}) that this implies 
$$\displaystyle \lim_{x\to 0}\myfrac{w(x)}{\phi_\gm^\Gw(x)}=k^*,
$$
and
$$\myint{\Gw}{}\left(w\CL^*_\gm\gz +\gz(g(u^*)-g(v))\right)d\gg_\gm^\Gw=k^*c_\mu\gz(0)\qquad\text{for all }\,\gz\in\BBX_\gm(\Gw).
$$
Since $u^*=w+v$ and
$$\myint{\Gw}{}\left(v\CL^*_\gm\gz +\gz g(v)\right)d\gg_\gm^\Gw=\myint{\Gw}{}\gz d(\gg_\gm^\Gw\gn^*)
+\myint{\prt\Gw}{}\gz d(\gb_\gm^\Gw\gl^*)\quad\text{for all }\,\gz\in\BBX_\gm(\Gw),
$$
 we infer
 \bel{3-D7}
  \BA {lll}
\myint{\Gw}{}\left(u^*\CL_\gm^*\gz+g(u^*)\gz\right)d\gg_\gm^\Gw=\myint{\Gw}{}\gz d(\gg_\gm^\Gw\gn^*)
+\myint{\prt\Gw}{}\gz d(\gb_\gm^\Gw\gl^*)+k^*c_\mu\gz(0)\quad\text{for all }\,\gz\in\BBX_\gm(\Gw).
\EA\ee
The last assertion is obvious.
\qeda\medskip

\begin{definition}\label{good} The triplet of measures $(\gn^*,\gl^*,k^*\gd_0)$ is called the reduced triplet  associated to $(\gn,\gl,k\gd_0)$. If $(\gn^*,\gl^*,k^*\gd_0)=(\gn,\gl,k\gd_0)$ the triplet is called $g$-good.
\end{definition}

\blemma{P2} Let $(\gn,\gl,k)$ and $(\gn',\gl',k')$ in $\frak M_+(\Gw;\gg^\Gw_\gm)\ti\frak M_+(\prt\Gw; \gb^\Gw_\gm)\ti\BBR_+$. If $\gn'\leq\gl$, $\gl'\leq\gl$ and $k'\leq k$ and $(\gn,\gl,k)=(\gn^*,\gl^*,k^*)$, then 
$(\gn',\gl',k')=(\gn'^*,\gl'^*,k'^*)$.
\es
{\bf Proof. } For $\ell>0$, let $u_\ell=u_{\ell,\gn,\gl,k}$  be the solution of (\ref{3-DX2}). We define similarly $u'_\ell=u'_{\ell,\gn',\gl',k'}$.  Then $u'_{\ell}\leq u_\ell$ for any $\ell>0$. Then $u_\ell\downarrow u^*$ and $u'_\ell\downarrow u'^*$  as $\ell\to\infty$ where $u^*\, u'^*$ are the solution of (\ref{1-B1}) with sources
$(\gn^*,\gl^*,k^*),\, (\gn'^*,\gl'^*,k'^*)$ respectively,  and these convergences hold in $L^1(\Gw,\gr^{-1}d\gg_\gm^\Gw)$ by the previous proposition.  Since $(\gn,\gl,k)=(\gn^*,\gl^*,k^*)$, then
$$\CL_\gm(u_\ell-u^*)+g_\ell(u_\ell)-g_\ell(u^*)= g(u^*)-g_\ell(u^*),
$$
and we deduce from \rprop{Kato} that
$$\myint{\Gw}{}(u_\ell-u^*)\gr^{-1}d\gg_\gm^\Gw+\myint{\Gw}{}|g_\ell(u_\ell)-g_\ell(u^*)|\eta d\gg_\gm^\Gw
\leq \myint{\Gw}{}(g(u)-g_\ell(u))\eta d\gg_\gm^\Gw.
$$
Because $|g_\ell(u_\ell)-g_\ell(u^*)|\leq |g_\ell(u_\ell)-g(u^*)|+g(u^*)-g_\ell (u^*)$ we get
$$\myint{\Gw}{}|g_\ell(u_\ell)-g_\ell(u^*)|\eta d\gg_\gm^\Gw\leq 2\myint{\Gw}{}(g(u^*)-g_\ell(u^*))\eta d\gg_\gm^\Gw\to 0\quad\text{as }\;\ell\to\infty.
$$
Since $g_\ell(u'_\ell)\leq g_\ell(u_\ell)$, it follows by Vitali's theorem that $g_\ell(u'_\ell)$ converges to $g(u'^*)$ in $L^1(\Gw,d\gg_\gm^\Gw)$. Letting $\ell\to\infty$ in the weak formulation of the equation satisfied by $u'_\ell$ we conclude that $u'^*$ verifies
$$
\myint{\Gw}{}\left(u'^*\CL_\gm^*\gz+g(u'^*)\gz\right)d\gg_\gm^\Gw=\myint{\Gw}{}\gz d(\gg_\gm^\Gw\gn')
+\myint{\prt\Gw}{}\gz d(\gb_\gm^\Gw\gl')+k'c_\mu\gz(0)\quad\text{for all }\,\gz\in\BBX_\gm(\Gw).
$$
This implies the claim.\qeda\medskip

As a consequence we have

\bprop{opt} The triplet $(\gn^*,\gl^*,k^*\gd_0)$ is the largest $g$-good triplet smaller than  $(\gn,\gl,k\gd_0)$. 
\es

\blemma{SUP} Let $(\gn,\gl,k)$ in $\frak M_+(\Gw;\gg^\Gw_\gm)\ti\frak M_+(\prt\Gw; \gb^\Gw_\gm)\ti\BBR_+$. The two next statements are equivalent:\\
(i) The triplet $(\gn,\gl,k)$ is $g$-good.\\
(ii) For any $\ge>0$, $0\leq k'\leq k$, $(\gn_\ge,\gl_\ge,k')$ is $g$-good.
\es
{\bf Proof. }  We recall that $\gn_\ge=\chi_{_{\Gw_\ge}}\gn$ and $\gl_\ge=\chi_{_{\prt\Gw\cap B^c_\ge}}\gl$. \\
(i) implies (ii) by \rlemma{P2}. \\
Conversely, if $(\gn_\ge,\gl_\ge,k')$ is $g$-good for any $\ge>0$ and $k'\in [0,k]$, let $u_{\ge,k'}$ be the solution of 
  \bel{3-D8}\left.\BA {lll}
\CL_\mu u+g(u)=\gn_\ge&\quad\text{in }\;\Gw\\[1mm]
\phantom{\CL_\mu+g(u)}u=\gl_\ge+k'\gd_0&\quad\text{on } \prt\Gw.
\EA\right.\ee
Then map $ (\ge,k')\mapsto u_{\ge,k'}$ is nonincreasing in $\ge$ and nondecreasing in $k'$.  There holds
 $$
\myint{\Gw}{}\left(u_{\ge,k'}\CL_\gm^*\gz+g(u_{\ge,k'})\gz\right)d\gg_\gm^\Gw=\myint{\Gw}{}\gz d(\gg_\gm^\Gw\gn_\ge)
+\myint{\prt\Gw}{}\gz d(\gb_\gm^\Gw\gl_\ge)+k'c_\mu\gz(0)\quad\text{for all }\,\gz\in\BBX_\gm(\Gw).
$$
From (\ref{2-C11}) we have that 
$$
\myint{\Gw}{}\left(u_{\ge,k'}\gr^{-1}+g(u_{\ge,k'})\eta\right)d\gg_\gm^\Gw\leq \myint{\Gw}{}\eta d(\gg_\gm^\Gw\gn_\ge)
+\myint{\prt\Gw}{}\eta d(\gb_\gm^\Gw\gl_\ge)+k'c_\mu \sup_{\Omega}\eta.
$$
Put $\displaystyle u=\lim_{(\ge,k')\to (0,k)}u_{\ge,k'}$. By the monotone convergence theorem, 
$$
\myint{\Gw}{}\left(u\gr^{-1}+g(u)\eta\right)d\gg_\gm^\Gw=\myint{\Gw}{}\eta d(\gg_\gm^\Gw\gn)
+\myint{\prt\Gw}{}\eta d(\gb_\gm^\Gw\gl)+kc_\mu\eta(0).
$$
Therefore $u_{\ge,k'}\to u$ in $L^1(\Gw,\gr^{-1}d\gg_\gm^\Gw)$ and $g(u_{\ge,k'})\to g(u)$ in $L^1(\Gw,d\gg_\gm^\Gw)$ as $\ge\to0^+$. Going to the limit in (\ref{3-D8}) yields the claim.\qeda\medskip

\nind \Remark The previous result is a particular case of the following result: If $\{(\gn_n,\gl_n,k_n)\}
\subset \frak M_+(\Gw;\gg^\Gw_\gm)\ti\frak M_+(\prt\Gw; \gb^\Gw_\gm)\ti\BBR_+$ is an increasing sequence of 
$g$-good triplet converging to $(\gn,\gl,k)\in \frak M_+(\Gw;\gg^\Gw_\gm)\ti\frak M_+(\prt\Gw; \gb^\Gw_\gm)\ti\BBR_+$, then $(\gn,\gl,k)$ is $g$-good.\medskip

\subsection{Capacitary framework, good measures and removable sets} 
In the sequel, we set $g(r)=g_p(r):=|r|^{p-1}r$ with $p>1$. The following a priori estimate of Keller-Osserman type is by now standard (see e.g. \cite {GuVe}, \cite{MaVe0}). 
\blemma{KE} Let $p>1$, $\gm\in\BBR$, $G\subset\BBR^N$ be a domain such that $0\notin G$. There exist constants $A>0$, $B\geq 0$ depending on $N$, $p$, $\gm$ such that any compact subset  $F$ of $\prt G$, possibly empty, and any  solution $v$ of 
  \bel{3-D8x}
  \left.\BA {lll}
\CL_\mu v+g_p(v)=0&\quad\text{in }\;G\\[1mm]
\phantom{\CL_\mu+g_p(v)}v=0&\quad\text{on }\;\prt G\setminus (\{0\}\cup F),\\
\EA\right.\ee
there holds
  \bel{3-D8y}
  \BA {lll}
|v(x)|\leq A\max\left\{|x|^{-\frac{2}{p-1}},\left(\dist (x,F)^{-\frac{2}{p-1}}\right)\right\}+B\quad\text{for all }\;x\in G.
\EA\ee
\es

\nind {\bf Proof of Theorem D.}    Since $g_p$ satisfies the uniform $\Gd_2$-condition, i.e. $K(|r|)$ is constant in inequality (\ref{1-B10}), if $(\gn,0,0)$, $(0,\gl,0)$ and $(0,0,k\gd_0)$
are $g_p$-good, then $(\gn,\gl,k\gd_0)$ is also $g_p$-good, and conversely. 
Assume now that $(\gn,\gl,0)$ is $g_p$-good, or, equivalently, for any $\ge>0$, $(\gn_\ge,\gl_\ge,0)$, is $g_p$-good. Let $u_\ge$ be the solution of (\ref{3-D8}) with $k'=0$. Let $\tilde \Gw_\ge$ be a smooth domain such that $\Gw_\ge\subset \tilde \Gw_\ge\subset \Gw_{\frac \ge 2}$. Then $\tilde u_\ge:=u_\ge\lfloor_{\tilde \Gw_\ge}$ satisfies 
 \bel{3-D9*}
  \left.\BA {lll}
\CL_\mu\tilde u_\ge+g_p(\tilde u_\ge)=\gn_\ge&\quad\text{in }\;\tilde\Gw_\ge\\[1mm]
\phantom{\CL_\mu+g_p(\tilde u_\ge)}\tilde u_\ge=\gl_\ge&\quad\text{on }\;\prt\Gw\cap\prt\tilde \Gw_\ge\\[1mm]
\phantom{\CL_\mu+g_p(\tilde u_\ge)}\tilde u_\ge=u_\ge&\quad\text{on }\;\prt\tilde \Gw_\ge\cap\Gw.
\EA\right.\ee
Furthermore $\frac{\gm}{|x|^2}$ is bounded in $\tilde \Gw_\ge$. Hence the Green operator $\BBG^{-\Gd +\gm|.|^{-2}}$ relative to $\tilde \Gw_\ge$ is equivalent of the one relative to $-\Gd$ and  $\gn_\ge\in \mathfrak M_+(\Gw;\gr)$.
Let $\tilde \Gw_{\ge,t}=\{x\in \tilde \Gw_\ge:\gr(x)>t\}$ and $\gn_{\ge,t}=\chi_{_{\tilde \Gw_{\ge,t}}}\gn_\ge$.
The bounded measure $\gn_{\ge,t}$ is $g_p$-good in $\tilde \Gw_\ge$. From \cite{BaPi}, this holds if and only if for any Borel set $K\subset\tilde \Gw_\ge$, 
$$Cap^{\BBR^N}_{2,p'}(K)=0\Longrightarrow \gn_{\ge,t}(K)=0.
$$
Assume now $E\subset\Gw$ is a compact set such that $Cap^{\BBR^N}_{2,p'}(E)=0$. 
Then $Cap^{\BBR^N}_{2,p'}(E\cap \tilde \Gw_{\ge,t})=0$ and thus $\gn_{\ge,t}(E\cap \tilde \Gw_{\ge,t})=0$. By the monotone convergence theorem, it implies
$$\displaystyle
\lim_{\ge\to 0}\lim_{t\to 0}\gn_{\ge,t}(E\cap \tilde \Gw_{\ge,t})=\lim_{\ge\to 0}\gn_\ge(E\cap \tilde \Gw_{\ge})
=\gn(E)=0.
$$
Similarly, using Marcus-V\'eron results on the boundary trace (see e.g. \cite{MaVe0}) $\gl$ is 
$g_p$-good if and only if $\gl_\ge$ vanishes on compact sets $E\subset\prt\tilde \Gw_\ge$ such that $Cap^{\BBR^{N-1}}_{\frac 2p,p'}(E)=0$. Clearly $\gl$ shares this property. \smallskip

Conversely, if $\gn$ (resp. $\gl$) vanishes on compact sets $E\subset\Gw$ (resp. $E\subset\prt\Gw$) such that 
$Cap^{\BBR^N}_{2,p'}(E)=0$ (resp. $Cap^{\BBR^N}_{\frac 2p,p'}(E)=0$), then $\gn_+$ (resp. $\gl_+$) has the same property. Hence we can assume that $\gn$ (resp. $\gl$) is nonnegative. Clearly $\gn_\ge$ (resp. $\gl_\ge$) shares also this property. If $0<t<\ge$ we denote by $\Gw_t^*$ a smooth domain such that $\Gw_\ge\subset \Gw_t^*\subset\Gw$ and 
$\Gw_t^*\cap B_{\frac t2}=0$ there exists an increasing sequence $\{\gn_{\ge,n}\}$ (resp. $\{\gl_{\ge,n}\}$) of positive bounded measures belong to $W^{-2,p}(\Gw)$ (resp. $W^{-\frac2p,p}(\prt\Gw)$) converging to $\gn_\ge$ (resp. ${\gl_\ge}$). 
The measures $\gn_{\ge,n}$ (resp. $\gl_{\ge,n}$) are $g_p$-good relatively to the open set $\Gw_t^*$. 
Therefore there exists a sequence of solutions $\{\tilde u_{\ge,t\,n}\}$ satisfying weakly
 \bel{3-D10}
  \left.\BA {lll}
\CL_\mu\tilde u_{\ge,t\,n}+g_p(\tilde u_{\ge,t\,n})=\gn_{\ge,n}&\quad\text{in }\;\Gw_t^*\\[1mm]
\phantom{\ \CL_\mu\tilde  +g_p(\tilde u_{\ge,t\,n})}\tilde u_{\ge,t\,n}=\gl_{\ge,n}&\quad\text{on }\;\prt\Gw\cap\prt\Gw_t^*\\[1mm]
\phantom{\CL_\mu +g_p(\tilde u_{\ge,t\,n})} \tilde u_{\ge,t\,n}=0&\quad\text{on }\;\prt\Gw_t^*\cap\Gw.
\EA\right.\ee
Letting $n\to\infty$, we infer that $\tilde u_{\ge,\gd\,n}$ increases and converges to the solution $\tilde u_{\ge,\gd}$ of 
  \bel{3-D11}
  \left.\BA {lll}
\CL_\mu\tilde u_{\ge,t}+g_p(\tilde u_{\ge,t})=\gn_{\ge}&\quad\text{in }\;\Gw_t^*\\[1mm]
\phantom{\CL_\mu+g_p(\tilde u_{\ge,\gd})}\tilde u_{\ge,t}=\gl_{\ge}&\quad\text{on }\;\prt\Gw\cap\prt\Gw_t^*\\[1mm]
\phantom{\CL_\mu+g_p(\tilde u_{\ge,\gd})}\tilde u_{\ge,t}=0&\quad\text{on }\;\prt\Gw_t^*\cap\Gw.
\EA\right.\ee
For $0<t<t'$, $\tilde u_{\ge,t}\geq \tilde u_{\ge,t'}$, hence $\tilde u_{\ge}:=\displaystyle\lim_{t\to 0^+}\tilde u_{\ge,t}$ satisfies 
 \bel{3-D12}
\myint{\Gw}{}\left(\tilde u_{\ge}\CL_\gm^*\gz+g_p(\tilde u_{\ge})\gz\right)d\gg_\gm^\Gw=\myint{\Gw}{}\gz d(\gg_\gm^\Gw\gn_\ge)
+\myint{\prt\Gw}{}\gz d(\gb_\gm^\Gw\gl_\ge),
\ee
for all $\gz\in\BBX_\gm^\Gw$ which vanishes in a neigborhood of $0$. We end the proof as in Theorem B. We first obtain that $\tilde u_{\ge}$ satisfies (\ref{3-D12}) for all $\gz\in\BBX_\gm^\Gw$, and then we let $\ge\to 0$ and conclude that $u:=\displaystyle\lim_{\ge\to 0}\tilde u_{\ge}$ satisfies
\bel{3-D13}
\left.\BA {lll}
\CL_\mu u+g_p(u)=\gn&\quad\text{in }\;\,\Gw\\[1mm]
\phantom{\CL_\mu+g_p(u)}u=\gl&\quad\text{on }\;\prt\Gw,
\EA
\right.\ee
hence $(\gn,\gl)$ is $g_p$ good.
\qeda\medskip

\nind{\bf Proof of Theorem E. }  A particular case of {\it Theorem E} that we will prove in {\it Theorem J} is that $0$ is a non-removable singularity if and only if $1<p<p^*_\gm$ for any $\gm\geq\gm_1$ and $N>2$, or $p>p^{**}_\gm$ with $N\geq 3$ and $\gm<1-N$.  

\nind (i) Assume $K\subset\Gw$ is compact. It follows from \cite[Theorem 3.1]{BaPi} that $Cap^{\BBR^N}_{2,p'}(K)=0$ is a necessary and sufficient condition for $K$ to be removable for the operator $\CL_\gm$ (and $p\geq \frac{N}{N-2}$ otherwise $K$ is empty).\smallskip

\nind (ii) Let $K\subset\prt\Gw\setminus\{0\}$ be compact and, for $\ge>0$, $K_\ge=\{x\in\Gw:\dist(x,K)<\ge\}$. Assume $u$ is a function belonging to $L^1(\Gw\setminus K_\ge,\gr^{-1}d\gg_\gm^\Gw)\cap L^p(\Gw\setminus K_\ge,d\gg_\gm^\Gw)$ for any $\ge>0$ satisfying 
 \bel{3-D14}\BA {lll}
\myint{\Gw}{}\left(u\CL^*_\gm\gz +g_p(u)\gz\right)d\gg_\gm(x)=0,
\EA\ee
for any $\gz\in\BBX_\gm(\Gw)$ vanishing in a neighborhood of $K$. Taking a test function $\gz\in C^{2}(\overline\Gw)$ vanishing on $\prt\Gw$ and in a neighborhood of $K$ we infer by standard regularity theory that 
$u\in C^2(\overline\Gw\setminus(K\cup\{0\})$ is a strong solution of $\CL_\gm u+g_p(u)=0$ in $\Gw$ which vanishes on
$\prt \Gw\setminus(K\cup\{0\})$. Let $G\subset\Gw$ be a smooth domain such that $K$ is interior to $\prt G\cap\prt\Gw$ relatively to the induced topology on $\prt\Gw$ and such that $0\notin\overline G$. Then $\gm|x|^{-2}$ is bounded in $\overline G$. Then there exists $a>0$ and $b\in\BBR$ such that 
$$g_p(u)+\gm|x|^{-2} u\geq au_+^p-b.
$$
Set $m=\max\{u_+(x):x\in \prt G\cap\Gw\}$. Then implies that $v=\left(u-m-\left(\frac{b_+}a\right)^{\frac 1p}\right)$ satisfies $-\Gd v+av^p\leq 0$ in $G$ and vanishes on $\prt G\setminus K$. Since $Cap^{\BBR^{N-1}}_{\frac 2p,p'}(K)=0$
(and $p\geq \frac{N+1}{N-1}$ otherwise $K$ is empty), $v=0$ by \cite[Theorem 3.3]{MaVeX}, which implies 
$u\leq m+\left(\frac{b_+}a\right)^{\frac 1p}$ in $G$. Similarly $u$ is bounded from below in $G$ and it follows that 
(\ref{3-D14}) holds for all $\gz\in\BBX_\gm(\Gw)$. Hence $u=0$ by uniqueness.

Conversely, if $Cap^{\BBR^{N-1}}_{\frac 2p,p'}(K)>0$, then there exists a capacitary measure $\gl_K$ belonging to $W^{-\frac 2p,p}(\prt\Gw)$ with support in $K$ (see \cite[Chapter 1]{AH}). Since $\gl_K$ vanishes on Borel set with $Cap^{\BBR^{N-1}}_{\frac 2p,p'}$-capacity $0$, it is $g_p$-good and there exists a solution $u$ to 
 \bel{3-D15}\left.\BA {lll}
\CL_\gm u+g_p(u)=0\quad&\text{in }\;\Gw\\[1mm]
\phantom{\CL_\gm +g_p(u)}
u=\gl_K\quad&\text{on }\;\prt\Gw.
\EA\right.\ee
Hence $u$ satisfies (\ref{3-D14}) for all $\gz\in\BBX_\gm(\Gw)$ vanishing in a neighborhood of $K$. Hence $K$ is not removable. \smallskip

\nind (iii) If $K\subset \overline \Gw\setminus\{0\}$ is such that $Cap^{\BBR^{N-1}}_{\frac 2p,p'}(K\cap\prt\Gw)>0$ then 
$K\cap\prt\Gw$ is not removable by (ii). If $c^{\BBR^{N}}_{ 2,p'}(K\cap\Gw)>0$, then there exists an increasing sequence of compact sets $K_n\subset K\cap\Gw$ such that $c^{\BBR^{N}}_{ 2,p'}(K_n)>0$. Hence $K_n$ is not removable, and clearly $K$ inherits the same property as it contains $K_n$. \smallskip

\nind (iv) If $0\in K\subset\prt\Gw$ and $K\setminus\{0\}\neq \emptyset$ and assume that any solution of (\ref{1-B13}) is identically $0$, in particular any solution which vanishes on $\prt\Gw\setminus\{0\}$ is zero. By {\it Theorem J} this is ensured only if $p\geq p^*_\gm$. If $Cap^{\BBR^{N-1}}_{\frac 2p,p'}(K)>0$, then either $p<\frac{N+1}{N-1}$, thus $K\setminus\{0\}$ contains at least one point which is not removable, or  $p\geq\frac{N+1}{N-1}$, and since $Cap^{\BBR^{N-1}}_{\frac 2p,p'}(K\setminus\{0\})=Cap^{\BBR^{N-1}}_{\frac 2p,p'}(K)>0>0$, there exists a compact subset $K'\subset K\setminus\{0\}$ such that $Cap^{\BBR^{N-1}}_{\frac 2,p'}(K')>0$. Hence $K'$, and therefore $K$, is not removable. This implies that if $K$ is removable one must have $p\geq p^*_\gm$ and $Cap^{\BBR^{N-1}}_{\frac 2p,p'}(K)=0$. 

  Conversely, if $p\geq p^*_\gm$, we will see at Theorem J that there exists no nonzero solution $u\in C(\overline\Gw\setminus\{0\})$ of $\CL_\gm u+g_p(u)=0$ vanishing on  $\prt \Gw\setminus\{0\}$. For $0<t<\ge$ we set 
  $K_t=\{x\in \Gw:\dist (x,K)<t\}$, $K_{\gd,\ge}=K_t\cap B_{2\ge}^c$ and 
$\Gw_{t,\ge}=\Gw\setminus \overline K_{t,\ge}$. We denote by $v_{t,\ge}$ the maximal solution of $\CL_\gm u+g_p(u)=0$ in $\Gw_{t,\ge}$ which vanishes on $\prt\Gw\setminus K_{t,\ge}$; hence it blows-up on $\prt K_{t,\ge}$ and it can be easily constructed by \rlemma {KE} by approximation with solutions with finite boundary value on $\prt K_{t,\ge}$. We also denote by $w_\ge$ the maximal solution of  the same equation in $\Gw_\ge:=\Gw\cap \overline B_\ge^c$ which vanishes on $\prt\Gw\setminus B_\ge$. It blows up on $\prt B_\ge\cap \Gw$. If $u$ is a solution of (\ref{1-B13}), it is dominated in $\Gw\setminus (\overline K_{t,\ge}\cup \overline B_\ge)$ by the supersolution 
$v_{t,\ge}+w_\ge$. When $t\to 0$, $v_{t,\ge}$ converges to the function $v_{0,\ge}$ which satisfies 
the equation in $\Gw_\ge$ and vanishes on $\prt\Gw_\ge\setminus K$. Since $c_{\frac{2}{p},p}^{\BBR^{N-1}}(K)=0$, there holds $c_{\frac{2}{p},p}^{\BBR^{N-1}}(K\cap B_{2\ge}^c)=0$. Therefore $v_{0,\ge}=0$. When $\ge\to 0$, 
$w_\ge$ decreases and converges to a solution of the equation in $\Gw$ which vanishes on $\prt\Gw\setminus\{0\}$, hence this limit is zero and consequently $u=0$.\smallskip

\nind (v) If $0\in K\subset\overline\Gw$ and $K\setminus\{0\}\neq \emptyset$ and any solution of (\ref{1-B13}) is identically $0$. Then $p\geq p^*_\gm$ as in (iv). Since $K\cap\Gw\neq \emptyset$ then any point in $K\cap\Gw$ is a removable singularity, hence $p\geq\frac{N}{N-2}$ (which implies $p>\frac{N+1}{N-1}$).
If $c_{2,p'}^{\BBR^N}(K\cap\Gw)>0$, there exists a compact set $K'\subset K\cap\Gw$ such that $c_{2,p'}^{\BBR^N}(K')>0$. Then $K'$ is not removable by {\it Theorem D}, hence $K$ is not removable too. If $c_{\frac 2p,p'}^{\BBR^{N-1}}(K\cap\prt\Gw)>0$, then $K$ is not removable as in (iv). 

Conversely assume that $p\geq p^*_\gm$, $c_{2,p'}^{\BBR^N}(K\cap\Gw)=0$, $c_{\frac 2p,p'}^{\BBR^{N-1}}(K\cap\prt\Gw)=0$ and $u$ satisfies  (\ref{1-B13}). For $0<t<\ge$,  we define $K_t=\{x\in \Gw:\dist (x,K\cap\prt\Gw)<t\}$, $K_{t,\ge}=K_\gd\cap B_{2\ge}^c$ as in (iv) and $\tilde K_{t,\ge}=\{x\in\Gw:\dist (x,K\cap\Gw)<t\}\cap \{x\in\Gw:\dist (x,\prt\Gw)>2\ge\}$. The functions $v_{t,\ge}$ and $w_\ge$ are defined as in (iv). We also denote 
by $\tilde v_{t,\ge}$ the maximal solution of $\CL_\gm+g_p(u)=0$ in $\Gw\setminus \tilde K_{t,\ge}$ which vanishes on $\prt\Gw$. Then $u\leq v_{t,\ge}+\tilde v_{t,\ge}+w_\ge$ in $\Gw\setminus (\overline K_{t,\ge}\cup \tilde K_{t,\ge}\cup \overline B_\ge$. When $t\to 0$, $v_{t,\ge}\to 0$ since $c_{\frac 2p,p'}^{\BBR^{N-1}}(K\cap B_\ge^c\cap\prt\Gw)=0$ and $\tilde v_{t,\ge}\to 0$ since $c_{2,p'}^{\BBR^N}(K\cap\{x\in\Gw:\dist (x,\prt\Gw)>2\ge\})=0$. Hence $u\leq w_\ge$ and we conclude as in (iv) by letting $\ge\to 0$.\qeda

\mysection{Isolated boundary singularities}
The study of boundary isolated singularities is based upon a technical framework which has been introduced by \cite{GmVe} in the case $\gm=0$. For the sake of completeness we recall this formalism, adapting it to our framework. Up to a rotation we assume that the inward normal direction to $\prt\Gw$ at $0$ is ${\bf e}_N=(0',1)\in\R^{N-1}\times\R$ and that the tangent hyperplane to $\prt\Gw$ at $0$ is $\prt\BBR^N_+=\BBR^{N-1}$. For $R>0$ set $B'_R=\{x'\in\BBR^{N-1}:|x'|<R\}$ and $D_R=B'_R\ti (-R,R)$. Then there exist $R>0$ and a $C^2$ function $\gth: B'_R\mapsto\BBR$ such that $\prt\Gw\cap D_R=\{x=(x',x_N): \, x_N=\gth(x')\ \text{for }x'\in B'_R\}$ and $\Gw\cap D_R=\{x=(x',x_N): \gth(x')<x_N<R\}$. Furthermore $\nabla \gth(0)=0$. Define the function $\Gth=(\Gth_1,...,\Gth_N)$ on $D_R$ by $y_j=\Gth_j(x)=x_j$ if $1\leq j\leq N-1$
and $y_N=\Gth_N(x)=x_N-\gth(x')$. Since $D\Gth(0)=Id$ we can assume that $\Gth$ is a diffeomorphism from $D_R$ onto $\Gth(D_R)$. Let $z$ be the harmonic extension of $h$ in $B_R\cap\Gw$ vanishing on $\Gw\cap \prt B_R$ and set
\bel{4-F1}
\BA {lll}\displaystyle
u(x)-z(x)=\tilde u(y),\; z(x)=\tilde z(y)\qquad\text{for all }x\in D^+_R=B'_R\ti[0,R).
\EA\ee
Denote by $(r,\gs)\in (0,\tilde r)\ti S^{N-1}$ the spherical coordinates in $\BBR^{N}$ and set 
\bel{4-F2}
\BA {lll}\displaystyle
\tilde u(y)=\tilde u(r,\gs)=r^{-a}v(t,\gs)\,,\;\tilde z(y)=\tilde z(r,\gs)=r^{-a}\CZ(t,\gs)\,,\; t=\ln r\,,\; a=\frac{2}{p-1}.
\EA\ee
Then $v$ is bounded and satisfies the following asymptotically autonomous equation in $(-\infty,r_0]\ti \BBS^{N-1}_+$
\bel{4-F3}\BA {lll}\displaystyle
(1+\ge_1(t,\cdot))v_{tt}+\left(N-2-2a+\ge_2(t,\cdot)\right) v_t+\left(a(a+2-N)-\gm+\ge_3(t,\cdot)\right)v+\Gd' v\\[2mm]
\phantom{}
+\langle \nabla' v,\ge_4(t,\cdot)\rangle+\langle \nabla' v_t,\ge_5(t,\cdot)\rangle+\langle \nabla'(\langle\nabla' v,{\bf e}_N\rangle),\ge_6(t,\cdot)\rangle+\gm\CZ-|v+\CZ|^{p-1}(v+\CZ)
=0,
\EA\ee
where $\Gd'$ is the Laplace-Beltrami operator on $\BBS^{N-1}$ and the $\ge_j$ satisfy the estimates
\bel{4-F4}
\BA {lll}\displaystyle
|\ge_j(t,\cdot)| +|\prt_t\ge_j(t,\cdot)|+|\nabla'\ge_{j}(t,\cdot)|\leq c_{17}e^t.
\EA\ee
As for $\CZ$ it verifies
\bel{4-F4a}
\BA {lll}\displaystyle
|\CZ(t,\cdot)| +|\prt_t\CZ(t,\cdot)|+|\nabla'\CZ(t,\cdot)|\leq c_{17}e^{at}.
\EA\ee
This is due to the fact that $|\gth(x')|=O(|x'|^2)$ near $0$. Furthermore, standard elliptic equations theory implies that
there holds, if $k+\ell\leq 3$,
\bel{4-F5}
\BA {lll}\displaystyle
\left|\frac{\prt^{k}\nabla'^{\ell}v}{\prt t^k}(t,\cdot)\right|\leq c_{18}\quad\text{in }(-\infty,r_0]\ti \BBS^{N-1}_+.
\EA\ee
\medskip

\nind{\bf Proof of Theorem F}. We denote by $\CS_{\gm,p}$ the set of functions satisfying
  \bel{4-F9}\left.\BA {lll}
-\Gd'\gw+\left(a(N-2-a)+\gm\right)\gw+g_p(\gw)=0&\quad\text{in }\ \BBS^{N-1}_+\\[2mm]
\phantom{-\Gd'+\left(a(N-2-a)+\gm\right)\gw+g_p(\gw)}
\gw=0&\quad\text{on }\prt \BBS^{N-1}_+,
\EA\right.\ee
where $a=\frac{2}{p-1}$. \\
(i) If $\gw$ is a solution it satisfies
\bel{4-F10}
  \BA {lll}
0=\myint{\BBS^{N-1}_+}{}\left(|\nabla'\gw|^2+\left(a(N-2-a)+\gm\right)\gw^2+|\gw|^{p+1}\right)dS\\
[4mm]\phantom{0}
\geq \myint{\BBS^{N-1}_+}{}\left(N-1+\left(a(N-2-a)+\gm\right)\gw^2+|\gw|^{p+1}\right)dS.
\EA\ee
If $N-1+a(N-2-a)+\gm\geq 0$, then necessarily $\gw=0$.  Next we have the following equivalences
\bel{4-F11}
  \BA {lll}
\phantom{NNNN}N-1+a(N-2-a)+\gm\geq0\Longleftrightarrow -\ga_+\leq a\leq-\ga_-\\
\phantom{N-1+a(N-2-a)+NNNN\gm\geq0}
\Longleftrightarrow\\[2mm]
\phantom{N}
\left\{\BA {lll}\text{(i)\; either }\;p\geq1-\frac{2}{\ga_-}= p^*_\gm,\\[2mm]
\!\text{(ii)\; or }\;1<p\leq 1-\frac{2}{\ga_+}=p^{**}_\gm\quad\text{provided }\,N\geq 3\text{ and }\;\gm_1\leq \gm< 1-N.
\EA\right.
\EA\ee
(ii) By minimization $\CS_{\gm,p}$ is not empty if the conditions (i) or (ii) of Theorem F are fulfilled, in which case 
$\CS_{\gm,p}$ has a unique positive element (see \cite{GmVe} for a similar situation). This unique positive element is denoted $\gw_\gm$.\\
(iii) The last statement follows an idea introduced in \cite{Ve0}. The hupper hemisphere admits the following representation
$$\BBS_+^{N-1}=\left\{x=((\sin\gf)\gs',\cos\gf):\gs'\in \BBS^{N-2}, \gf\in \left(0,\frac{\gp}{2}\right)\right\}.
$$
The surface measure $dS$ on $\BBS^{N-1}$ can be decomposed as 
$$dS(\gs)=(\sin\gf)^{N-2}dS'(\gs') d\gf 
$$
where $dS'$ is the surface measure on $\BBS^{N-2}$. If $h(\gs)=h(\gs',\gf)$ is defined on $\BBS^{N-1}$, we put
$$\overline h'(\gf)=\myfrac{1}{|\BBS^{N-2}|}\myint{\BBS^{N-2}}{}h(\gs',\gf)dS'(\gs').
$$
Let $\gw$ be an element of  $\CS_{\gm,p}$, then, by averaging (\ref{4-F9}), 
$$\BA{lll}
\myint{\BBS^{N-1}_+}{}\left(-\Gd'(\gw-\overline\gw')+\left(a(N-2-a)+\gm\right)(\gw-\overline\gw')+\left(g_p(\gw)-\overline{g_p(\gw)}'\right)\right) (\gw-\overline\gw')dS=0.
\EA$$
By monotonicity
$$\BA {lll}\myint{\BBS^{N-1}_+}{}\left(g_p(\gw)-\overline{g_p(\gw)}'\right)(\gw-\overline\gw')dS
=\myint{\BBS^{N-1}_+}{}\left(g_p(\gw)-g_p(\overline\gw')\right)(\gw-\overline\gw')dS\\[4mm]
\phantom{\myint{\BBS^{N-1}_+}{}\left(g_p(\gw)-\overline{g_p(\gw)}'\right)(\gw-\overline\gw')dS}
\geq 2^{1-p}\myint{\BBS^{N-1}_+}{}|\gw-\overline\gw'|^{p+1}dS.
\EA$$
The function $\gw-\overline\gw'$ is orthogonal to the first eigenspace of $-\Gd'$ in $H^{1,2}_0(\BBS^{N-1}_+)$. Since the second eigenvalue of $-\Gd'$ in $H^{1,2}_0(\BBS^{N-1}_+)$ in $2N$, we have
$$-\myint{\BBS^{N-1}_+}{}(\gw-\overline\gw')\Gd'(\gw-\overline\gw')dS\geq 2N\myint{\BBS^{N-1}_+}{}(\gw-\overline\gw')^2dS.
$$
Hence
$$\myint{\BBS^{N-1}_+}{}\left(\left(a(N-2-a)+\gm+2N\right)(\gw-\overline\gw')^2+2^{1-p^{\phantom{p^p}}}\!\!\!\!|\gw-\overline\gw'|^{p+1}\right)dS\leq 0.
$$
Hence, if $a(N-2-a)+\gm+2N\geq 0$ it follows that $\gw-\overline\gw'=0$. The polynomial 
$$P_2(X):=X^2+(2-N)X-\gm-2N
$$
admits two real roots provided $\gm\geq-(\frac{N+2}{2})^2:=\gm_2$,  which are expressed by
$$a_-=\frac N2 -1+\sqrt{\gm-\gm_2},\quad\text{ }\; a_+=\frac N2 -1-\sqrt{\gm-\gm_2},
$$
and 
$$P_2\left(\frac{2}{p-1}\right)\leq 0 \Longleftrightarrow a_+\leq \myfrac{2}{p-1}\leq a_-.
$$
Note that $a_+a_->0$ if and only if $-2N>\gm$. Furthermore $P_2(-\ga_-)<0$ and $P_2(-\ga_+)<0$. Then there holds
$$\BA {lll}
(i)\; \,\,\text{if } \gm\geq1-N\;&\text{then }\; a_+<-\ga_+<0<-\ga_-<a_-\Longrightarrow\tilde p^*_\gm<p^*_\gm,\\[2mm]
(ii) \,\,\,\text{if }N\geq 3\ \&\ -2N\leq \gm< 1-N\;&\text{then }\; a_+<0<-\ga_+\leq-\ga_-<a_-\Longrightarrow \tilde p^*_\gm<p^*_\gm<p^{**}_\gm, \\[2mm]
(iii)\,\,\text{if }N\geq 9\ \&\ \gm_1\leq \gm< -2N\;&\text{then }\; 0<a_+<-\ga_+\leq-\ga_-<a_-\Longrightarrow \tilde p^*_\gm<p^*_\gm<p^{**}_\gm<\tilde p^{**}_\gm,
\EA $$
where, we recall it,
$$p^*_\gm=1-\frac{2}{\ga_-},\quad\tilde p^*_\gm=1+\frac{2}{a_-},\quad\tilde p^{**}_\gm=1+\frac{2}{a_+}, \quad p^{**}_\gm=1-\frac{2}{\ga_+}.
$$
 Therefore $\gw-\overline\gw'=0$ if the following conditions are satisfied
 \bel{4-F12}
  \BA {lll}
(i)\quad\text{when } \,\gm\geq 1-N\;\text{ and }\quad \tilde p^*_\gm\leq p< p^*_\gm,\quad\quad\quad\quad\quad
\quad\quad\quad\\[2mm]
(ii)\quad\text{when } N\geq 3\,,\;-2N\leq \gm< 1-N\,\text{ and either}\, \tilde p^*_\gm\leq p< p^*_\gm \text{ or }\,
p^{**}_\gm<p,\\[2mm]
(iii)\quad \text{when } N\geq 9\,\text{ and }\;\gm_1\leq \gm< -2N\text{ and either}\, \tilde p^*_\gm\leq p< p^*_\gm \text{ or }p^{**}_\gm<p\leq \tilde p^{**}_\gm.
\EA\ee
If one of the above conditions is fulfilled,   
$\gw$ depends only on the variable $\gf\in [0,\frac\gp2]$. It satisfies
 \bel{4-F13}
  \left.\BA {lll}
-\myfrac{1}{\sin^{N-2}\gf}\left(\gw_\gf\sin^{N-2}\gf\right)_\gf+\left(a(N-2-a)+\gm\right)\gw+g_p(\gw)=0\quad\text{ in }\,(0,\frac\gp2)\\[3mm]\phantom{------------------\ }
\gw_\gf(0)=0\,,\;\gw(\frac\gp2)=0.\quad
\EA\right.\ee
Define the operator 
$$\CB(\psi)=-\myfrac{1}{\sin^{N-2}\gf}\left(\sin^{N-2}\gf\psi_\gf\right)_\gf,$$
among functions $\psi$ in the space $\CH_\CB\subset  C^2([0,\frac\gp2])$ satisfying $\psi_\gf(0)=0$ and $\psi(\frac\gp2)=0$. 
The first eigenvalue of $\CB$ in $\CH_\CB$ is $N-1$ and the second in $2N$. Since $g_p$ is nonnecreasing, it is known (see e.g. \cite{Ber}) that the constant sign solutions $\gw_p$ and $-\gw_p$ lie on a branch of bifurcation issued from $N-1$ and there exists no other bifurcation when the parameter $a(N-2-a)+\gm$ belongs to $(N-1,2N]$. 
This implies $\CS_{\gm,p}=\{\gw_p,-\gw_p,0\}$ and ends the proof. \qeda
\medskip 

For proving {\it Theorems G, H, I, J} we recall here the following technical results \cite[Theorem 5.1]{GmVe} related to the solutions of (\ref{1-B14}) satisfying
\bel{4-Fws}
\lim_{x\to 0}|x|^{\frac{2}{p-1}}u(x)=0.
 \ee
 The statement is easily adapted from the one of the above mentioned theorem. We denote by $\gl_k=\{k(k+N-2:k\in\BBN^*\}$ the set of eigenvalues of $-\Gd'$ in $H^{1,0}(\BBS^{N-1}_+)$. Any separable  $\CL_\gm$-harmonic function in 
 $\BBR^N_+$ vanishing on $\prt\BBR^N_+\setminus\{0\}$ endows the form
 \bel{4-Fws1}
x\mapsto u(x)=u(r,\gs)=r^{\ga_k}\phi_k(\gs)\qquad (r,\gs)\in \BBR_+\ti \BBS^{N-1}_+,
 \ee
 where $\phi_k\in ker(\Gd'+\gl_k I)$ and $\ga_k$$=\ga_{k-}$ or $\ga_{k+}$ the smallest and the largest root of 
 \bel{4-Fws2}
\ga^2+(N-2)\ga-\gl_k-\gm=0,
 \ee
which exist for some $k\geq 1$ if and only if $\gm\geq\gm_k:=\gm_1+N-1-\gl_k$. Note that $\ga_{k-}\leq 0$ for all  $k\in\BBN^*$ and $\ga_{k+}\leq 0$ if and only if $\gm\geq -\gl_k$ (which imposes $N\geq8k(k+\sqrt{2k(k-1)})$).
 
 \bth{T51} Assume $\gm\geq\gm_1$, $1<p<p^*_\gm$ and $h\in C^3(\prt\Gw)$. If $u\in C(\overline\Gw\setminus\{0\})\cap C^2(\Gw)$ is a solution of (\ref{1-B14})  satisfying (\ref{4-Fws}) and
\smallskip
  
  \nind (A)   either  $u_-(x)=O(|x|^{-\frac{2}{p-1}+\gd})$ near $x=0$, for some $\gd>0$, \smallskip
  
  \nind (B)  or $N=2$ and $\Gw$ is locally a straight line near $x=0$, \smallskip
  
  \nind (B)  or $-\frac{2}{p-1}$ is not equal to some  $\ga_{k-}$ for some $k\in\BBN^*$.\smallskip
  
  \nind Then \smallskip
  
  \nind (i) either $u$ is the weak solution of (\ref{1-B14h}),\\
  \nind (ii) or there exist an integer $k\in\BBN^*$ such that $-\ga_{k-}<\frac{2}{p-1}$  and a nonzero spherical harmonic $\psi_k$ of degree $k$ such that 
    \bel{4-Fws3}\displaystyle
\lim_{x\to 0}r^{\ga_{k-}}\tilde u(r,\gs)=\psi_k(\gs).
 \ee
  \es

\subsection{Proof of Theorems G, H, I and J}
 Because of (\ref{4-F5}) the negative trajectory of $v$ in $C_0^1(\overline{\BBS^{N-1}_+})$ which is defined by
\bel{4-F6*}
 \CT_-(v)=\bigcup_{t\leq r_0-1}\{v(t,.)\},
\ee
is relatively compact in the $C^2(\overline{\BBS^{N-1}_+})$-topology. The limit set $\CE_v$ of $ \CT_-(v)$ at $-\infty$ defined by
\bel{4-F6}
\CE_v=\bigcap_{\gt\leq r_0-1}\overline{\bigcup_{t\leq \gt}\{v(t,.)\}}^{C_0^1(\overline{\BBS^{N-1}_+})},
 \ee
is non-empty. Since $1<p<p^*_\gm$ and  $\gm\geq \gm_1$, there holds
\bel{4-F7}
p<\frac{N+2}{N-2}.
 \ee
Thus the coefficient of $v_t$ in (\ref{4-F3}) is not zero (asymptotically, when $t\to-\infty$). Then energy damping holds and, in the same way as in \cite{GmVe} up to a shift of $\gm$ in the coefficient of $v$ in (\ref{4-F3}), we obtain
 \bel{4-F8}
\myint{-\infty}{r_0-1}\myint{\BBS^{N-1}_+}{}v_t^2d\gs dt<\infty. 
\ee
 Combining this estimate with (\ref{4-F5}) and some standard manipulations (see \cite{GmVe}) implies that 
 $$\norm{v_t(t,.)}_{C^1(\overline{\BBS^{N-1}_+})}+\norm{v_{tt}(t,.)}_{C(\overline{\BBS^{N-1}_+})}\to 0\quad\text{as }t\to-\infty.$$
 Hence $\CE_v$ is a compact connected component of $\CS_{\gm,p}$. \medskip
 
\nind{\bf Proof of Theorem G. } 
  If $\gw$ is nonnegative, either $\CE_v=\{\gw_p\}$ and (\ref{1-B19}) holds or 
   \bel{4-F9*}\displaystyle
\lim_{t\to-\infty} \norm{v(t,.)}_{C^2(\overline{\BBS^{N-1}_+})}\to 0\quad\text{as }t\to-\infty.
\ee
 If this holds, it follows by \rth{T51}-A that either $u=0$ or (\ref{4-Fws3}) is verified for some $k\geq 1$. Since any spherical harmonics of degree at least two changes sign $k$ must be equal to $1$. Then 
 $$\tilde u(x)=\ell\phi_\gm(x)(1+o(1))\quad\text{as }x\to 0,
 $$
which is  (\ref{1-B20}).\qeda
\bcor{A} Let $\gm_1\leq \gm$ and $1<p<p_\gm^*$. Then for any $h\in C^3(\Gw)$, $h\geq 0$ there exists only one 
solution of (\ref{1-B14}) with a strong singularity at $x=0$, that is satisfying (\ref{1-B19}).
\es
{\bf Proof.} It is a consequence of {\it Theorem G} that the limit of $u_{\ell\gd_0,h}$ of the solution of (\ref{1-B14hk}) when $\ell\to\infty$ is a solution which satisfies (\ref{1-B19}). The method of proof of uniqueness is due to Marcus and V\'eron \cite{MaVeZ}. The minimal solution of (\ref{1-B14hk}) with a strong singularity at $x=0$ is defined by
$$\underline u_{\infty,h}:=\lim_{\ell\to\infty}u_{\ell\gd_0,h}.$$
For constructing the maximal solution we define the sequence $\overline u_{n,h}$  of solutions of 
  \bel{4-U1}
   \left.\BA {lll}\displaystyle
\CL_\gm \overline u_{n,h}+g_p(\overline u_{n,h})=0&\quad\text{in }\,\, \Gw\cap \overline B^c_{\frac 1n}\\[2mm]
\phantom{\CL_\gm +g_p(\overline u_{n,h})}\overline u_{n,h}=h&\quad\text{on }\,\prt\Gw\cap  B^c_{\frac 1n}\\
\phantom{\CL_\gm +g_p(\overline u_{n,h})}\overline u_{n,h}=cn^{\frac{2}{p-1}}&\quad\text{on }\,\Gw\setminus \prt B_{\frac 1n},
\EA \right.\ee
where $c>0$ is some constant large enough. Then $u_{\ell\gd_0,h}\leq \overline u_{n,h}$. By convexity there holds
$$\overline u_{n,h}-u_{\ell\gd_0,h}\leq \overline u_{n,0}-u_{\ell\gd_0,0}.
$$
By monotonicity $\{\overline u_{n,h}\}$ decreases and converges to the maximum solution $\overline u_{\infty,h}$ of
(\ref{1-B19}) and there holds
$$\overline u_{\infty,h}-\underline u_{\infty,h}\leq \overline u_{\infty,0}-\underline u_{\infty,0}.
$$
Furthermore, by (\ref{1-B19}), there exists $K=K(p,\gm,\Gw)>1$ such that 
  \bel{4-U2}
   \BA {lll}\displaystyle
\overline u_{\infty,0}\leq K\underline u_{\infty,0}.
\EA \ee 
 If we assume that $\overline u_{\infty,0}>\underline u_{\infty,0}$, then,  again by convexity, the function 
 $$U=\overline u_{\infty,0}-\frac{1}{K}\left(\underline u_{\infty,0}-\overline u_{\infty,0}\right)
 $$
 is a supersolution for problem (\ref{1-B14}) smaller than $\underline u_{\infty,0}$. The function 
  $$U^*=\left(\frac{1}{2K}+\frac{1}{2}\right)\overline u_{\infty}
  $$
  is a supersolution of the same problem (\ref{1-B14})  smaller than $U$. By a classical result valid for a wide class of quasilinear equations there exists $V$ solution of the problem such that $U^*\leq V\leq U$. In particular $V$ has a strong blow-up at $x=0$ and it is smaller than the minimal solution $\underline u_{\infty}$, contradiction.\qeda

 \medskip
 
\nind{\bf Proof of Theorem H. }\  Since $\CE_v$ is a connected subset of the discrete set $\CS_{\gm,p}$ which has three connected components $(\{\gw_p\},\{-\gw_p\},\{0\})$ by Theorem F-(1) either ({\ref{1-B21}) or (\ref{4-Fws}) holds. Since $p>\tilde p^*_\gm $, $-\frac{2}{p-1}$ which necessarily larger $\ga_{1-}$ satisfies either  $\frac{2}{p-1}<\ga_{2-}$ or, if  $\frac{2}{p-1}>\ga_{2+}$ in the case $N\geq 9$ and $\gm<-2N$ and $\frac{2}{p-1}$ is not equal to any $\ga_{k-}$ or $\ga_{k+}$ for $k>2$ by the equation. Hence, by \rth{T51}, (\ref{1-B20}) holds.\qeda
 \medskip
 
 \nind\Remark If $p=\tilde p^*_\gm $ or $p=\tilde p^{**}_\gm $ the method shows that either ({\ref{1-B21}) or (\ref{4-Fws}) holds. Since it is the spectral case always difficult to handle we cannot prove that (\ref{1-B20}) also holds,  a fact that we conjecture.\smallskip
 
\smallskip
 
\nind{\bf Proof of Theorem I. } 
 The two statements obey a totaly different approach. \smallskip

\nind Statement 1- is a consequence of the theory of analytic functionals developped by in \cite{Si1}, \cite{Si2} and applied to Emden-Fowler equations in \cite{B-VV}. The key point is to consider the equation (\ref{4-F3}) satisfied by $v(t,.)=r^{-\frac{2}{p-1}}\tilde u(r,\gs)$ in $(-\infty,r_0)\ti \BBS^{N-1}_+$ and to verify that, as a function of $v$, it is {\it real analytic}. Hence $p$ must be an odd integer. If $1<p<p^*_\gm$ the only possibility is $p=3$ which is in the range if $N+4\sqrt{\gm-\gm_1}$. If $\gm<1-N$ and $p>p^{**}_\gm$ there are infinitely many possibilities for $p$.\smallskip

\nind Statement 2- The convergence to one element of $\CS_{\gm,p}$ follows from the fact that this set of solutions of (\ref{1-B24}) is discrete. If $\prt\Gw$ is locally a close graph near $0$ 
the paper \cite{CMV} which use Sturmian arguments and the Jordan closed curve theorem.  If $\gw=0$, as quoted in \cite[Theorem 5.1-C2]{GmVe} we perform a reflexion through $\prt\Gw$ near $0$ and apply the result of \cite[Lemma 2.1]{CMV}, the shift of the coefficient by $\gm$ playing no role.\qeda
  \medskip

\nind{\bf Proof of Theorem J. }\ 
   Since $p>1$, any solution $u$ of (\ref{1-B14}) satisfies the estimate of \rlemma {KE} under the following form
  \bel{4-J1}\displaystyle
|u(x)|\leq A|x|^{-\frac{2}{p-1}}\qquad\text{for all }x\in\overline\Gw\setminus\{0\}. 
 \ee
 If $p>p_\gm^*$, then $-\frac{2}{p-1}>\ga_-$, therefore $u(x)=o(\phi_\gm(x))$ near $x=0$. Let $u_+$ be the solution of 
  \bel{1-B14h+}\left.\BA {lll}
\CL_\gm u+g_p(u)=0&\quad\text{in }\;\Gw\\[1mm]
\phantom{\CL_\gm +g_p(u)}
u=h_+&\quad\text{on }\;\prt\Gw.
\EA\right.\ee
For any $\ge>0$,  $u_++\ge\phi_\gm$ is a supersolution of $\CL_\gm u+g_p(u)=0$, larger than $u$ near $x=0$.
 Then that $u\leq u_++\ge\phi_\gm$ and, letting $\ge\to 0$ then $u\leq u_+$. Similarly $u$ is larger than $-u_--\ge\phi_\gm$, where $u_-$ is the solution of (\ref{1-B14h+}) with $h_+$ replaced by $h_-$. Letting $\ge \to 0$ yields 
 $-u_-\leq u\leq u_+$. It follows by the method of Theorem B that $u$ is the weak solution of (\ref{1-B14h}).\smallskip
 
  \nind If $p=p_\gm^*$ and $\gm=\gm_1$, then similarly $u(x)=o(\phi_\gm(x))$ near $x=0$ and the result follows by the same method. \smallskip 
  
   \nind Finally, if $p=p_\gm^*$ and $\gm>\gm_1$, then (\ref{4-F7}) holds. Using the variable 
 $t=\ln r$ and $v(t,.)=r^{\frac2{p-1}}\tilde u(r,.)$ we obtain from the previous energy method that 
 $$\CE_v\subset\CS_{\gm,p}=\{0\}.
 $$
 Hence $u$ satisfies (\ref{4-Fws}). Since $p=p_\gm^*$, $\frac{2}{p-1}=-\ga_-$. Hence $u=o(\phi_\gm)$ near $0$ and the conclusion follows as in the previous cases. \qeda

 \medskip
 
 \noindent{\bf Acknowledgements} H. Chen is is supported by NNSF of China, No: 12071189 and 12001252, by the
Jiangxi Provincial Natural Science Foundation, No: 20202BAB201005, 20202ACBL201001.

\end{document}